# DISTRIBUTED NONPARAMETRIC FUNCTION ESTIMATION: OPTIMAL RATE OF CONVERGENCE AND COST OF ADAPTATION[*]

By T. Tony Cai and Hongji Wei

*University of Pennsylvania*

Distributed minimax estimation and distributed adaptive estimation under communication constraints for Gaussian sequence model and white noise model are studied. The minimax rate of convergence for distributed estimation over a given Besov class, which serves as a benchmark for the cost of adaptation, is established. We then quantify the exact communication cost for adaptation and construct an optimally adaptive procedure for distributed estimation over a range of Besov classes.

The results demonstrate significant differences between nonparametric function estimation in the distributed setting and the conventional centralized setting. For global estimation, adaptation in general cannot be achieved for free in the distributed setting. The new technical tools to obtain the exact characterization for the cost of adaptation can be of independent interest.

**1. Introduction.** Distributed statistical estimation and inference are becoming increasingly important as in many applications data can be necessarily distributed at different locations due to the size constraint or privacy and security concerns. Such a setting arises in a range of medical, financial, and business applications. With distributed data, separate statistical analyses need to be carried out at individual sites and then the results are transmitted to and aggregated at a central location in order to make the final statistical decision. For large-scale data analysis, communication costs can be expensive and become the main bottleneck in statistical practice. It is important to understand the interplay between communication constraints and statistical accuracy, as well as how to design optimal estimation and inference procedures under communication constraints.

There has been an increasing amount of recent literature on distributed estimation when the communication budget is limited. For example, Zhang et al. (2013); Garg et al. (2014); Braverman et al. (2016); Han et al. (2018); Zhu and Lafferty (2018); Szabó and van Zanten (2018); Barnes et al. (2019); Cai and Wei (2020a); Szabo and van Zanten (2020) considered information-theoretical limits under communication constraints for various distributed estimation problems, such as Gaussian mean estimation, linear regression and nonparametric

[*]The research was supported in part by NSF grants DMS-1712735 and DMS-2015259 and NIH grants R01-GM129781 and R01-GM123056.

*MSC 2010 subject classifications:* Primary 62F30; secondary 62B10, 62F12

regression. Optimality results have been established under different communication constraints. Besides theoretical analysis, progress has also been made on developing practical methodologies for distributed estimation. See, for example, Kleiner et al. (2014); Deisenroth and Ng (2015); Lee et al. (2017); Diakonikolas et al. (2017); Jordan et al. (2019); Battey et al. (2018); Fan et al. (2019). Further literature review is given in Section 1.4.

In this paper we study distributed minimax and distributed adaptive nonparametric estimation under communication constraints in a decision theoretical framework. In the conventional non-distributed settings, adaptation has been a central goal for nonparametric function estimation. it is well-known that adaptive estimation can be achieved for free under a range of global losses such as the integrated squared error over a wide collection of Besov classes (Donoho and Johnstone, 1995; Johnstone, 2017). Indeed, it is possible to adaptively achieve superefficiency for free (Cai, 2008). However, in the distributed settings, adaptation becomes more difficult and involved due to the additional communication constraints. A rate-optimal adaptive algorithm needs to perform well statistically while efficiently compressing the information from the local machines to the central learner. Intuitively, the difficulty arises from the fact that only limited amount of information can be transmitted and information that is critical for estimation over one function class might not be essential for estimation over another. In such a setting, it is easy to imagine that achieving adaptation over a collection of function classes requires more communication budget than what is needed for a given function class in the minimax setting.

The primary goal of the present paper is to precisely characterize the communication cost of adaptation for distributed nonparametric function estimation. We first establish the minimax rate of convergence for distributed estimation over a given Besov class, which serves as a benchmark for the cost of adaptation when the smoothness parameters are unknown. We then quantify the exact cost of adaptation and construct an optimally adaptive procedure for distributed nonparametric estimation over a range of Besov classes.

1.1. *Distributed estimation framework.* We begin by introducing a general framework for distributed estimation by giving a formal definition of transcript, distributed estimator, and independent distributed protocol. Let $\mathcal{P} = \{P_\theta : \theta \in \Theta\}$ be a parametric family of distributions supported on space $\mathcal{X}$, where $\theta \in \Theta$ is the parameter of interest. Suppose there are $m$ local machines and a central machine, where the local machines contain the observations and each local machine has access only to data in that machine, and the central machine produces the final estimator of $\theta$ under the communication constraints between the local and central machines. More precisely, suppose we observe i.i.d. random samples drawn from a distribution $P_\theta \in \mathcal{P}$:

$$X_i \stackrel{\text{iid}}{\sim} P_\theta, \quad i = 1, \ldots, m,$$

where the $i$-th local machine has access to $X_i$ only.

On each machine, because of limited communication budget, the observation $X_i$ on the $i$-th local machine needs to be processed to a uniquely decodable string $Z_i$ by a (possibly random) function $Z_i : \mathcal{X} \to \bigcup_{b=1}^{\infty} \{0,1\}^b$. The resulting string $Z_i = Z_i(X_i)$, which is called the **transcript** from the $i$-th machine, is then transmitted to the central machine. Here we denote the length of transcript $Z_i$ as $|Z_i|_l$, which indicates the communication cost for sending this transcript. Finally, a **distributed estimator** $\hat{\theta}$ is constructed on the central machine based on the transcripts $Z_1, Z_2, ..., Z_m$,

$$\hat{\theta} = \hat{\theta}(Z_1, Z_2, ..., Z_m).$$

The above scheme to obtain a distributed estimator $\hat{\theta}$ is called an **independent distributed protocol**. Within an independent distributed protocol, the transcripts from each local machine only depend on its local observations and no information is exchanged between the local machines. There are also other types of distributed protocols with more interactive communication schemes in the literature (Zhang et al., 2013). In the present work we focus on independent distributed protocol. Define $L(\hat{\theta}) \triangleq \sum_{i=1}^{m} |Z_i|_l$ as the total communication cost for distributed estimator $\hat{\theta}$. The class of distributed protocols with total communication budgets $B$ can be defined as

$$\mathcal{A}_T(B) = \{(\hat{\theta}, Z_1, Z_2, ..., Z_m) : L(\hat{\theta}) \leq B\}.$$

The above classes of distributed protocol imposes uniform hard upper bounds on the length of transcripts, that is, the (total) length of transcripts are constrained to be less than a certain value given any possible observations. It is sometimes worthwhile to consider transcripts with variable length in order to gain possible adaptation to the data. In such settings, we introduce a class of distributed protocols with the expected total communication budgets for the family $\mathcal{P}$:

$$(1) \qquad \mathcal{A}_E(B, \Theta) = \{(\hat{\theta}, Z_1, Z_2, ..., Z_m) : \sup_{\theta \in \Theta} \mathbb{E}_{P_\theta} L(\hat{\theta}) \leq B\}$$

where the expected total communication cost is uniformly bounded by $B$ under any data generating distribution $P_\theta \in \mathcal{P}$.

As usual, the estimation accuracy of a distributed estimator $\hat{\theta}$ is measured by the mean squared error (MSE), $\mathbb{E}_{P_\theta} \|\hat{\theta} - \theta\|_2^2$, where the expectation is taken over the randomness in both the data and construction of the transcripts and estimator. As in the conventional decision theoretical framework, a quantity of particular interest in distributed learning is the minimax risk for the distributed protocols

$$\inf_{\hat{\theta} \in \mathcal{A}_E(B, \Theta)} \sup_{P_\theta \in \mathcal{P}} \mathbb{E}_{P_\theta} \|\hat{\theta} - \theta\|_2^2,$$

which characterizes the difficulty of the distributed learning problem under the expected total communication constraints $\mathcal{A}_E(B, \Theta)$. Similarly $\mathcal{A}_E(B, \Theta)$ can be replaced by other class of distributed protocols to illustrate minimax risk under other kind of communication constraints. In a rigorous decision theoretical formulation of distributed learning, the communication constraints are essential. Without the constraints, one can always output the original data from the local machines to the central machine and the problem is then reduced to the usual centralized setting.

1.2. *Distributed estimation.* We consider distributed minimax and adaptive estimation for the Gaussian sequence model and white noise model. For the white noise model, the goal is to recover the unknown function based on the noisy observations collected on $m$ machines, where on the $i$-th machine, for $1 \leq i \leq m$, one observes a Gaussian process,

$$dY_i(t) = f(t)dt + \frac{\epsilon}{\sqrt{n}}dW_i(t) \quad t \in [0,1], i = 1, 2, ..., m. \tag{2}$$

Here $\frac{\epsilon}{\sqrt{n}}$ is the noise level and $W_i(t), i = 1, 2, ..., m$ are independent standard Wiener process. The $i$-th machine has access to $Y_i(t)$ only. The goal is to recover the unknown function $f$ based on the distributed observed processes $Y_1(t), Y_2(t), ...., Y_m(t)$.

In the conventional centralized setting, wavelet methods (Donoho and Johnstone, 1994; Hall et al., 1999; Cai, 1999) have been shown to be a powerful tool for nonparametric function estimation as it decomposes a function into a structured wavelet series and a nonparametric estimation problem is then transformed into a Gaussian sequence estimation problem. Motivated by the equivalence between the white noise model and Gaussian sequence model, we begin by focusing on the following distributed Gaussian sequence estimation problem. Suppose there are $m$ machines, on $i$-th machine we have i.i.d Gaussian observations

$$X_{i,jk} = \theta_{jk} + \sigma z_{i,jk}, \quad j = 0, 1, 2, ...; k = 1, 2, ..., 2^j \tag{3}$$

where $z_{i,jk} \stackrel{\text{iid}}{\sim} N(0,1)$ for $i = 1, 2, ..., m; j = 0, 1, ...; k = 1, 2, ..., 2^j$, the noise level $\sigma$ known. The $i$-th machine can only access to $X_i \triangleq (X_{i,jk})_{j \geq 0, k=1,2,...,2^j}$ only. The goal is to estimate $\theta \triangleq (\theta_{i,jk})_{j \geq 0, k=1,2,...,2^j}$ under the mean-squared error

$$R(\hat{\theta}, \theta) = \|\hat{\theta} - \theta\|_2^2 = \sum_{j=0}^{\infty} \sum_{k=1}^{2^j} (\hat{\theta}_{jk} - \theta_{jk})^2.$$

We consider estimation over a collection of Besov classes $\mathcal{B}_{p,q}^{\alpha}(M)$ with $\alpha, p, q, M > 0$, where $\mathcal{B}_{p,q}^{\alpha}(M)$ is defined as the set of sequences $\theta$ satisfying $|\theta|_{b_{p,q}^{\alpha}} \leq M$ with the Besov sequence seminorm $|\theta|_{b_{p,q}^{\alpha}}$ given by

$$|\theta|_{b_{p,q}^{\alpha}} \triangleq \left( \sum_{j=0}^{\infty} \left( 2^{js} \left( \sum_{k=1}^{2^j} |\theta_{jk}|^p \right)^{1/p} \right)^q \right)^{1/q}. \tag{4}$$

Here $s = \alpha + 1/2 - 1/p > 0$ and $1 \leq p, q \leq \infty$, with the obvious replacement of the corresponding $\ell_p$ or $\ell_q$ norms to $\ell_\infty$ norms when $p, q = \infty$. The Besov sequence norm $|\theta|_{b_{p,q}^{\alpha}}$ is equivalent to the Besov function norm on the original function $f$; see, for example, Meyer (1992). Therefore, the distributed Gaussian sequence model (3) is equivalent to the white noise model (2). In the classical centralized setting, the Gaussian sequence model is also known to be a good proxy to study estimation of a function under the nonparametric regression model.

1.3. *Main contributions.* For estimation under the Gaussian sequence model (3) with communication constraints, a distributed estimation procedure, called seq-MODGAME, is proposed, and its rate of convergence under the communication constraints is derived. A matching lower bound is established to show that the seq-MODGAME procedure is optimal. The upper and lower bounds together yield the sharp optimal rate of convergence for estimation over a Besov class $\mathcal{B}_{p,q}^\alpha(M)$:

$$\mathcal{R}_E(B, \mathcal{B}_{p,q}^\alpha(M)) \triangleq \inf_{\hat{\theta} \in \mathcal{A}_E(B, \mathcal{B}_{p,q}^\alpha(M))} \sup_{\theta \in \mathcal{B}_{p,q}^\alpha(M)} \|\hat{\theta} - \theta\|_2^2.$$

where $\mathcal{A}_E(B, \mathcal{B}_{p,q}^\alpha(M))$ is the set of distributed protocols under the expected total communication constraints defined in (1). The same optimal rate holds for the white noise model. To the best of our knowledge, this is the first exact minimax rate of convergence for the distributed nonparametric function estimation. In comparison, the existing results have at least a logarithmic gap in the upper and lower bounds and are for more specialized parameter spaces such as a Hölder or Sobolev class.

We then quantify the exact communication cost for adaptation and construct an optimally adaptive procedure for distributed estimation over a range of Besov classes. Our analysis shows interesting phenomena. In the classical non-distributed setting, it is well known that adaptation can be achieved for free in terms of global risk measures such as the mean integrated squared error over a wide collection of Besov classes. See, for example, Donoho and Johnstone (1995); Johnstone (2017). However, in the distributed setting, our results show that there are unavoidable additional communication costs for any adaptive procedure over a collection of Besov classes. Specifically, the results provide a sharp characterization for the communication costs for adaptation, where it is shown that $O(m^3)$ total additional bits are necessary and sufficient to achieve the adaptation over a wide collection of Besov classes. In addition, a local thresholding procedure is constructed and is shown to be the most communication-efficient among all adaptive distributed estimators. Our newly proposed local thresholding procedure requires no prior knowledge on the range of the smoothness parameters, and is able to automatically achieve statistical adaptation over a wide collection of Besov classes $\mathcal{B}_{p,q}^\alpha(M)$ with $p \geq 2$ at the guaranteed minimum communication cost. The analysis on adaptive estimation makes significant improvement over existing results. The new technical tools used to obtain the exact characterization for the cost of adaptation can be of independent interest.

1.4. *Related literature.* Distributed nonparametric function estimation has been investigated in the recent literature. Zhu and Lafferty (2018) studied distributed minimax rate of convergence for the white noise model over the Sobolev classes with a logarithmic gap between the upper and lower bounds. Szabó and van Zanten (2018) derived distributed minimax rate for nonparametric regression under the integrated squared error and sup-norm error losses over the Hölder classes and Sobolev classes, also with a logarithmic gap between the upper and lower bounds. The paper also showed that adaptation is possible within the range $\alpha \in [\alpha_{\min}, \alpha_{\max})$ where $\alpha_{\min}$ depends on the given communication budget.

Szabo and van Zanten (2020) considered a two-point adaptation problem for distributed nonparametric estimation and showed that two-point adaptation is impossible when the smoothness indices of the two function classes are both larger than a certain threshold. It also proposed an adaptive distributed protocol that achieves statistical adaptation over a range of Sobolev classes with the smoothness indices below a certain threshold, while at the same time transmitting the minimal number of bits, up to a logarithmic gap. Szabo and van Zanten (2020) provided a clear solution when two-point adaptation can be achieved without additional communication cost. However, it is not clear whether adaptation is possible with additional communication budgets under the same settings. In comparison, we provide a more general lower bound for the communication cost for adaptive distributed estimators over a collection of Besov classes and construct an estimator that is adaptive over a wider range of parameter spaces at the guaranteed minimum communication cost.

1.5. *Organization of the paper.* We finish this section with notation, definitions, and some assumptions that will be used in the rest of the paper. Section 2 establishes the optimal rate of convergence for distributed Gaussian sequence estimation and Section 3 characterizes the communication cost of adaptation and introduces adaptive distributed procedures. The numerical performance of the proposed distributed estimators is investigated in Section 4 and further research directions are discussed in Section 5. For reasons of space, we only prove lower bounds for communication cost of adaptive estimators in Section 6 and defer the proofs of other main results and the technical lemmas to the supplementary material Cai and Wei (2020b).

1.6. *Notation, definitions, and assumptions.* For simplicity, in later sections we denote $n_j = 2^j$ be the number of coefficients at the $j$-th resolution level. For any positive integers $n, N$, let $[n] \triangleq \{1, 2, ..., n\}$ and $n \bmod N$ be the remainder of $n$ divided by $N$. For any $a \in \mathbb{R}$, let $\lfloor a \rfloor$ denote the floor function (the largest integer not larger than $a$). Unless otherwise stated, we shorthand $\log a$ as the base 2 logarithmic of $a$. For any $a, b \in \mathbb{R}$, let $a \wedge b \triangleq \min\{a, b\}$ and $a \vee b \triangleq \max\{a, b\}$. We use $a = O(b)$ or equivalently $b = \Omega(a)$ to denote there exist a constant $C > 0$ such that $a \leq Cb$, and we use $a \asymp b$ to denote $a = O(b)$ while $b = O(a)$. For any vector $a$, denote by $\|a\| \triangleq \sqrt{\sum_k \left(a^{(k)}\right)^2}$ its $l_2$ norm. For any finite set $S$, let $\text{card}(S)$ denote the cardinality of $S$. Define the density of a Gaussian distribution with mean 0 and standard deviation $\sigma$ as

$$\phi_\sigma(x) = \frac{1}{\sqrt{2\pi}\sigma} e^{-\frac{x^2}{2\sigma^2}}.$$

Throughout the paper, we shall assume $s = \alpha + 1/2 - 1/p > 0$. This condition is necessary for the estimation problem to be well-formulated. When $s \leq 0$, the closure of the Besov ball $\mathcal{B}^\alpha_{p,q}(M)$ is not compact and the compactness of the closure of the parameter space is a necessary condition for consistent estimation under the homoskedastic Gaussian sequence model. See Ibragimov and Khasminskii (1997) and Johnstone (2017, Theorem 5.7). Moreover, we assume $M \geq \sigma$. Otherwise estimation over the Besov ball $\mathcal{B}_{p,q}^\alpha(M)$ is trivial as the simple estimator $\hat{\theta} = 0$ is optimal.

**2. Minimax Optimal Rate of Convergence.** In this section we study the minimax rate of convergence for estimating the mean of a Gaussian sequence $\theta \in \mathcal{B}_{p,q}^\alpha(M)$ under the expected total communication constraint:

$$
(5) \qquad \inf_{\hat{\theta} \in \mathcal{A}_E(B, \mathcal{B}_{p,q}^\alpha(M))} \sup_{\theta \in \mathcal{B}_{p,q}^\alpha(M)} \mathbb{E}\|\hat{\theta} - \theta\|^2
$$

where we assume the parameters $\alpha, p, q, M$ are known in an oracle setting.

If there is no communication constraint, or equivalently we are in a centralized setting, Donoho and Johnstone (1998) pointed out the minimax rate of convergence over Besov classes is

$$
\inf_{\hat{\theta}} \sup_{\theta \in \mathcal{B}_{p,q}^\alpha(M)} \mathbb{E}\|\hat{\theta} - \theta\|^2 \asymp M^{\frac{2}{2\alpha+1}} \left(\frac{\sigma^2}{m}\right)^{\frac{2\alpha}{2\alpha+1}}.
$$

However, when the communication constraints take effect, there will be a loss of statistical accuracy thus the optimal rate of convergence (5) will further depend on the expected total communication cost $B$.

We first introduce a distributed estimation procedure satisfying the communication-constraint and provide an upper bound for its statistical performance. A matching lower bound on its minimax risk is then established. The upper and lower bounds together unveil a sharp minimax rate of convergence and the optimality of the proposed estimator.

2.1. *Optimal procedure.* We begin with the construction of an estimation procedure under the communication constraints and provide a theoretical analysis of the proposed procedure. The construction of the following procedure, called *seq-MODGAME*, is inspired by the MODGAME procedure proposed in Cai and Wei (2020a) for distributed Gaussian mean estimation. However, unlike the simple Gaussian mean estimation problem considered in Cai and Wei (2020a), the magnitude of each coordinate of $\theta$ is not known as a priori because within Besov space $\mathcal{B}_{p,q}^\alpha(M)$, the constraint on the Besov norm (4) is imposed on the whole vector, but not individual entries. Therefore, to estimate a mean vector $\theta \in \mathcal{B}_{p,q}^\alpha(M)$ under Gaussian sequence model (3), one needs a more sophisticated quantization strategy than the MODGAME procedure proposed in Cai and Wei (2020a).

We first define several useful functions and quantities. Define localization encoding function $g : \mathbb{Z} \to \bigcup_{k=1}^\infty \{0,1\}^k$ by the following rule:

- $g(0) = $ "0".

- When $x$ is a positive integer, let $k$ be the length of its binary representation, and define $g(x)$ to be a string starting with "1", followed by $k$ zeros and then followed by the binary representation of $x$. For example, $g(1) = $ "101" and $g(8) = $ "100001000".
- When $x$ is a negative integer, let $k$ be the length of the binary representation of $-x$, and define $g(x)$ to be a string starting with "11", followed by $k-1$ zeros and then followed by the binary representation of $-x$. For example, $g(-1) = $ "111" and $g(-8) = $ "110001000".

The function $g(x)$, as an encoding mechanism, has two main properties. First, it is a prefix code, thus uniquely decodable (Blahut and Blahut, 1987). We denote $g^{-1}$ as its inverse function (decoding function). Second, the length of $g(x)$ is guaranteed to be no larger than $2\log(|x|+1)+3$, which means that its length is adaptive to the magnitude of $x$. We will see that $g(x)$ plays an important role in the construction of the transcripts with variable length under the communication constraints.

As in the conventional centralized setting, we estimate the coordinates of the vector $\theta = (\theta_{j,k}) \in \mathcal{B}_{p,q}^\alpha(M)$ from its noisy observation up to a certain resolution level $j_{\max}$ and truncate all the components above $j_{\max}$ to zero. Note that when the communication budget is insufficient, the estimation accuracy in the distributed setting is not as good as in the centralized setting. So we first decide the maximal resolution level $j_{\max}$, and precision parameter $\delta$ according to communication budget $B$ and other model parameters. At those resolution levels lower than $j_{\max}$, we estimate each entry in an optimal way so that the stochastic error is roughly $O(\delta)$. At those higher resolution levels, we just truncate all entries to zero. The advanced communication strategy used in the procedure is the key to the optimality results.

We are now ready to introduce the seq-MODGAME procedure in detail. It is divided into two cases: $B < \left(\frac{M}{\sigma}\right)^{\frac{2}{2\alpha+1}}$ and $B \geq \left(\frac{M}{\sigma}\right)^{\frac{2}{2\alpha+1}}$.

**Case 1:** $B < \left(\frac{\Lambda_0 M}{\sigma}\right)^{\frac{2}{2\alpha+1}}$.

Let $\delta$ be a precision parameter calculated by

$$\delta \triangleq \Lambda_0 M B^{-(\alpha+1/2)}$$

where $\Lambda_0 > 0$ is a large tuning parameter. Let $j_{\max}$ be the maximal resolution level, defined as

$$j_{\max} \triangleq \max\left\{j : M \cdot 2^{-j(\alpha+1/2)} \geq \delta\right\}.$$

In this case, only one local machine is needed to sent transcripts to the central machine.

**First step: Generate the transcripts on the first local machine.** On the first local machine (who can access to data $X_1$), the output transcript $Z_1$ is the collection of the "crude localization" strings $Z_{1,jk}$, $0 \leq j \leq j_{\max}, k \in [n_j]$ where $Z_{1,jk}$ is defined as

$$Z_{1,jk} = g(\lfloor X_{1,jk}/\delta \rfloor).$$

**Second step: Generate the distributed estimator $\hat{\theta}^O$ on the central machine.**
The central machine can receive $Z_{1,jk}$, $0 \leq j \leq j_{\max}, k \in [n_j]$ from the first local machine. The final estimate $\hat{\theta}^O$ is given by

$$\hat{\theta}^O_{jk} = g^{-1}(Z_{1,jk}) \cdot \delta \text{ if } 0 \leq j \leq j_{\max}, k \in [n_j]$$

$$\hat{\theta}^O_{jk} = 0 \text{ if } j > j_{\max}, k \in [n_j]$$

**Case 2:** $B \geq \left(\frac{\Lambda_0 M}{\sigma}\right)^{\frac{2}{2\alpha+1}}$

Let $u$ be a parameter and $\delta$ be the precision parameter. They are calculated by

$$u \triangleq \left((\Lambda_0 M/\sigma)^{-\frac{1}{\alpha+1}} B^{\frac{2\alpha+1}{2\alpha+2}}\right) \wedge m, \quad \delta = \sigma/\sqrt{u},$$

and let $j_{\max}$ be the maximal resolution level, defined as

$$j_{\max} \triangleq \max\left\{j : M \cdot 2^{-j(\alpha+1/2)} \geq \delta\right\}.$$

In this case, with the help of communication strategy introduced in Cai and Wei (2020a), each entry of $\theta$ at lower resolution levels can be estimated in the most communication-efficient way so that their estimation errors is roughly of order $\delta$.

**First step: Generate the transcripts on the local machines.**

1. On the first machine (which has access to data $X_1$), the output transcript $Z_1$ is the collection of the "crude localization" strings $Z_{1,jk}$, $0 \leq j \leq j_{\max}, k \in [n_j]$ where $Z_{1,jk}$ is defined as

$$Z_{1,jk} = g(\lfloor X_{1,jk}/\sigma \rfloor);$$

2. On the $i$-th machine where $2 \leq i \leq 1 + \lfloor \log^2 u \rfloor$, the output transcript $Z_i$ is the collection of the "finer localization" strings $Z_{i,jk}$, $0 \leq j \leq j_{\max}, k \in [n_j]$ where $Z_{i,jk}$ is defined as

$$Z_{i,jk} = g(\lfloor X_{i,jk}/\sigma \rfloor \mod \lfloor \log u \rfloor);$$

3. On the $i$-th machine where $2 + \lfloor \log^2 u \rfloor \leq i \leq u$ the output transcript $Z_i$ is the collection of the 'refinement" strings $Z_{i,jk}$, $0 \leq j \leq j_{\max}, k \in [n_j]$ where $Z_{i,jk}$ is defined as

$$Z_{i,jk} = \lfloor X_{i,jk}/\sigma \rfloor \mod 8.$$

4. On the $i$-th machine where $u < i \leq m$, the local machine does not output anything.

**Second step: Generate the distributed estimator $\hat{\theta}$ on the central machine.** The central machine receives the transcripts $Z_1, Z_2, ..., Z_u$ from the local machines. Note that the code words in $Z_1, Z_2, ..., Z_u$ are all uniquely decodable, thus those transcripts can be decomposed into short strings $Z_{i,jk}$ for $i \in [u], j \in J_i, k \in [n_j]$.

The final estimator $\hat{\theta}^O$ is constructed as follows.

- For each $0 \leq j \leq j_{\max}$, $k \in [n_j]$:

  1. Because $g(x)$ is uniquely decodable, from $Z_{1,jk} = g(\lfloor X_{i,jk}/\sigma \rfloor)$ one can recover the value of $\lfloor X_{i,jk}/\sigma \rfloor$. Let $I_{jk}^a$ be an left-closed-right-open interval of length $u$ defined as

  $$I_{jk}^a \triangleq \begin{cases} \left[\lfloor X_{i,jk}/\sigma \rfloor - \frac{\lfloor \log u \rfloor - 1}{2}, \lfloor X_{i,jk}/\sigma \rfloor + \frac{\lfloor \log u \rfloor + 1}{2}\right) & \text{if } \lfloor \log u \rfloor \text{ is an odd number} \\ \left[\lfloor X_{i,jk}/\sigma \rfloor - \frac{\lfloor \log u \rfloor}{2}, \lfloor X_{i,jk}/\sigma \rfloor + \frac{\lfloor \log u \rfloor}{2}\right) & \text{if } \lfloor \log u \rfloor \text{ is an even number} \end{cases}$$

  2. Denote $z_{ik}^b \triangleq \arg\max_{z'} \sum_{i=2}^{\lfloor \log^2 u \rfloor + 1} \mathbb{I}_{\{Z_{i,jk} = z'\}}$ be the mode statistic among the $Z_{i,jk}, 2 \leq i \leq \lfloor \log^2 u \rfloor + 1$. Note that the length of $I_{jk}^a$ is $\lfloor \log u \rfloor$, so there will be exactly one integer $x_{jk}^b \in I_{jk}^a$ that satisfies

  $$x_{jk}^b \bmod \lfloor \log u \rfloor = g^{-1}(z_{ik}^b).$$

  Let $I_{jk}^b$ be an interval of length 3 defined by

  $$I_{jk}^b \triangleq [x_{jk}^b - 1, x_{jk}^b + 1].$$

  3. Let $p^h$ be the proportion of those refinement strings whose value is equal to $g(x_{jk}^b - 2 \bmod 8)$:

  $$p^h \triangleq \frac{1}{u - 1 - \lfloor \log^2 u \rfloor} \sum_{i=\lfloor \log^2 u \rfloor + 2}^{u} \mathbb{I}_{\{Z_{i,jk} = g(x_{jk}^b - 2 \bmod 8)\}}$$

  Define a function

  $$h_{jk}(y) \triangleq \sum_{l=-\infty}^{\infty} \int_{x_{jk}^b - 2 + 8l}^{x_{jk}^b - 1 + 8l} \phi_1(x - y) dx$$

  It is easy to see that $h_{jk}(y)$ is a strictly decreasing function on $I_{jk}^b$. Let $h_{jk}^{-1}(y)$ be the inverse function of $h_{jk}(y)$ which maps $h_{jk}(I_{jk}^b)$ to $I_{jk}^b$. The estimate is calculated by

  $$\hat{\theta}_{jk}^O = \begin{cases} (x_{jk}^b + 1)\sigma & \text{if } p^h \leq h_{jk}(x_{jk}^b + 1) \\ h_{jk}^{-1}(p^h)\sigma & \text{if } h_{jk}(x_{jk}^b + 1) < p^h < h_{jk}(x_{jk}^b - 1) \\ (x_{jk}^b - 1)\sigma & \text{if } p^h \geq h_{jk}(x_{jk}^b - 1) \end{cases}$$

- For each $j \geq j_{\max} + 1, k \in [n_j]$, set

$$\hat{\theta}_{jk}^O = 0.$$

The following theorem provides the theoretical guarantee for the communication cost of $\hat{\theta}^O$, as well as an upper bound for its statistical performance.

THEOREM 1. *If $\Lambda_0$ is set to be a sufficient large constant such that $\Lambda_0 > (24\alpha+64)^{\alpha+1/2}$, then the estimator $\hat{\theta}^O \in \mathcal{A}_E(B, \mathcal{B}_{p,q}^\alpha(M))$ and there exists a constant $C > 0$ such that*

$$
(6) \quad \sup_{\theta \in \mathcal{B}_{p,q}^\alpha(M)} \mathbb{E}\|\hat{\theta}^O - \theta\|^2 \leq C \cdot \begin{cases} M^2 B^{-2\alpha} & \text{if } B < \left(\frac{M}{\sigma}\right)^{\frac{2}{2\alpha+1}} \\ M^{\frac{2}{\alpha+1}} \left(\frac{\sigma^2}{B}\right)^{\frac{\alpha}{\alpha+1}} & \text{if } \left(\frac{M}{\sigma}\right)^{\frac{2}{2\alpha+1}} \leq B < \left(\frac{M}{\sigma}\right)^{\frac{2}{2\alpha+1}} m^{\frac{2\alpha+2}{2\alpha+1}} \\ M^{\frac{2}{2\alpha+1}} \left(\frac{\sigma^2}{m}\right)^{\frac{2\alpha}{2\alpha+1}} & \text{if } B \geq \left(\frac{M}{\sigma}\right)^{\frac{2}{2\alpha+1}} m^{\frac{2\alpha+2}{2\alpha+1}} \end{cases}
$$

*for all $2 \leq p \leq \infty, 0 < q \leq \infty, \alpha > 0, M > 0$*

REMARK 1. The proposed distributed estimator $\hat{\theta}^O$ satisfies expected total communication constraint, which is weaker than other types of constraint considered in the literature. The reason we work on this type of communication constraint is to illustrate the main idea and omit unnecessary complication when presenting the estimator. With suitable modification, the estimator can be made to satisfy other kinds of communication constraint, say, a fixed/hard total communication constraint or an equally assigned communication constraint on each single local machines.

For example, the following proposition provides a quick look on how $\hat{\theta}^O$ satisfies fixed/hard total communication constraint with high probability.

PROPOSITION 1. *With probability at least $1 - \exp(-B/18)$, we have*

$$L(\hat{\theta}^O) < 2B.$$

That is, the proposed estimator $\hat{\theta}^O$ satisfies the total communication constraint $2B$ with high probability. Note that the additional factor on the communication constraint doesn't affect the rate of convergence given in Theorem 1, therefore the estimator is still rate-optimal.

2.2. *Lower bound analysis.* Section 2.1 gives a detailed construction of the seq-MODGAME procedure for distributed Gaussian sequence estimation and provides a theoretical guarantee for the estimator in Theorem 1. In this section we shall show that the estimator $\hat{\theta}^O$ is indeed rate optimal among all estimators satisfying the total communication constraints by proving that the upper bound in Equation (6) cannot be improved. The following theorem gives a lower bound on the minimax risk under the expected total communication constraints.

THEOREM 2. *There exists a constant $c > 0$ such that*

$$(7) \quad \mathcal{R}_E(B, \mathcal{B}_{p,q}^\alpha(M)) \geq c \cdot \begin{cases} M^2 B^{-2\alpha} & \text{if } B < \left(\frac{M}{\sigma}\right)^{\frac{2}{2\alpha+1}} \\ M^{\frac{2}{\alpha+1}} \left(\frac{\sigma^2}{B}\right)^{\frac{\alpha}{\alpha+1}} & \text{if } \left(\frac{M}{\sigma}\right)^{\frac{2}{2\alpha+1}} \leq B < \left(\frac{M}{\sigma}\right)^{\frac{2}{2\alpha+1}} m^{\frac{2\alpha+2}{2\alpha+1}} \\ M^{\frac{2}{2\alpha+1}} \left(\frac{\sigma^2}{m}\right)^{\frac{2\alpha}{2\alpha+1}} & \text{if } B \geq \left(\frac{M}{\sigma}\right)^{\frac{2}{2\alpha+1}} m^{\frac{2\alpha+2}{2\alpha+1}} \end{cases}$$

*for all $0 < p \leq \infty, 0 < q \leq \infty, \alpha > 0, M > 0$*

The lower bound given in Theorem 2 is proved by constructing simultaneous tests $\theta_{jk} = 0$ vs $\theta_{jk} = \delta$ for all $j \leq J$, $k = 1, 2, ..., 2^j$, with pre-specified choices of $\delta$ and $J$. Then by strong data processing inequalities, we can prove that at least a proportion of entries cannot be accurately estimated. The detailed proof is deferred to the supplementary material Cai and Wei (2020b).

Theorems 1 and 2 together establish the minimax rate for distributed Gaussian sequence estimation:

$$(8) \quad \mathcal{R}_E(B, \mathcal{B}_{p,q}^\alpha(M)) \asymp \begin{cases} M^2 B^{-2\alpha} & \text{if } B < \left(\frac{M}{\sigma}\right)^{\frac{2}{2\alpha+1}} \\ M^{\frac{2}{\alpha+1}} \left(\frac{\sigma^2}{B}\right)^{\frac{\alpha}{\alpha+1}} & \text{if } \left(\frac{M}{\sigma}\right)^{\frac{2}{2\alpha+1}} \leq B < \left(\frac{M}{\sigma}\right)^{\frac{2}{2\alpha+1}} m^{\frac{2\alpha+2}{2\alpha+1}} \\ M^{\frac{2}{2\alpha+1}} \left(\frac{\sigma^2}{m}\right)^{\frac{2\alpha}{2\alpha+1}} & \text{if } B \geq \left(\frac{M}{\sigma}\right)^{\frac{2}{2\alpha+1}} m^{\frac{2\alpha+2}{2\alpha+1}} \end{cases}.$$

where $2 \leq p \leq \infty, q \leq \infty, \alpha > 0, M > 0$. The results also show that the distributed estimator $\hat{\theta}^O$ proposed in Section 2.1 is rate optimal under the total communication constraints.

The theorem also suggests that in order to achieve the centralized rate of convergence, which is of order $M^{\frac{2}{2\alpha+1}} \left(\frac{\sigma^2}{m}\right)^{\frac{2\alpha}{2\alpha+1}}$, a communication cost of order $\left(\frac{M}{\sigma}\right)^{\frac{2}{2\alpha+1}} m^{\frac{2\alpha+2}{2\alpha+1}}$ is sufficient and necessary.

REMARK 2. Similar as the optimal rate of convergence for distributed univariate Gaussian mean estimation Cai and Wei (2020a), the minimax rate (8) can be divided into three phases: localization ($B < \left(\frac{M}{\sigma}\right)^{\frac{2}{2\alpha+1}}$), refinement ($\left(\frac{M}{\sigma}\right)^{\frac{2}{2\alpha+1}} \leq B < \left(\frac{M}{\sigma}\right)^{\frac{2}{2\alpha+1}} m^{\frac{2\alpha+2}{2\alpha+1}}$), and optimal-rate ($B \geq \left(\frac{M}{\sigma}\right)^{\frac{2}{2\alpha+1}} m^{\frac{2\alpha+2}{2\alpha+1}}$). The minimax rate decreases quickly in the localization phase, when the communication constraints are extremely severe; then it decreases slower in the refinement phase, when there are more communication budgets; finally the minimax rate coincides with the centralized optimal rate (Donoho and Johnstone, 1998) and stays the same, when there are sufficient communication budgets. The value for each additional bit decreases as more bits are allowed.

REMARK 3. As mentioned in the introduction, distributed minimax estimation was considered in Zhu and Lafferty (2018) for the Hölder classes and in Szabo and van Zanten (2020) for the Sobolev classes. These two types of function classes are special cases of the Besov classes with the Hölder class being $B^\alpha_{\infty,\infty}$ and Sobolev class being $B^\alpha_{2,\infty}$. Furthermore, in both Zhu and Lafferty (2018) and Szabo and van Zanten (2020), the existing upper bound and lower bound are sub-optimal (with a poly-logarithmic gap to the optimal rate of convergence (8)). In contrast, the minimax rate given in (8) is sharp for a wide collection of Besov spaces.

**3. Adaptive Gaussian sequence estimation.** The minimax rate of convergence established in Section 2 provides an important benchmark for the evaluation of the performance of distributed Gaussian sequence estimators. However, the estimator $\hat{\theta}^O$, in spite of its statistical optimality and communication efficiency, requires explicit knowledge of the smoothness parameters which are typically unknown in practice. The optimal seq-MODGAME procedure proposed in Section 2 highly depends on the prior knowledge on the parameter space $\mathcal{B}^\alpha_{p,q}(M)$ so that local machines efficiently transmit useful information when the communication budget is limited. It is evident from the construction and theoretical analysis that the estimator $\hat{\theta}^O$ designed for one Besov class $B^\alpha_{p,q}(M)$ with a given smoothness parameter $\alpha$ would perform poorly over another Besov class $B^{\alpha'}_{p,q}(M)$ with a different smoothness parameter $\alpha'$. Therefore, the estimator $\hat{\theta}^O$ is not practical for real applications because the model parameters are typically unavailable.

This naturally leads to the important question of adaptive distributed estimation: Is it possible to construct a single distributed estimator, satisfying the communication constraints and not depending on the smoothness parameters, that achieves the optimal rate of convergence simultaneously over a wide collection of Besov classes $B^\alpha_{p,q}(M)$? In the conventional centralized setting, the answer is affirmative. That is, one can achieve adaptation for free for estimating a Gaussian sequence over a collection of Besov classes $B^\alpha_{p,q}(M)$ under the mean squared error.

Adaptive estimation in the centralized setting has been a major goal in the classical nonparametric function estimation literature. In particular, wavelet thresholding is well known to be a powerful technique to achieve adaptivity. For example, Donoho and Johnstone (1995); Abramovich et al. (2006) proposed adaptive term-by-term thresholding methods and Cai (1999); Cai and Zhou (2009) introduced data-driven block thresholding procedures to achieve optimal rate of convergence over a wide collection of Besov spaces. In contrast, little has been understood on how to construct a communication-efficient adaptive estimator for most distributed estimation problems, including but not limited to distributed Gaussian sequence estimation. It is interesting and practically important to investigate the interplay between communication constraints and adaptation for distributed estimation problems.

In this section we address the following questions: how to construct a data-driven distributed estimation procedure that can achieve the centralized optimal rate with communication cost as small as possible? Can adaptation be achieved for free? If not, what is the cost of adaptation?

It was shown in Section 2 that, for distributed estimation over the Besov class $\mathcal{B}_{p,q}^{\alpha}(M)$, one needs at least $\Omega\left(\left(\frac{M}{\sigma}\right)^{\frac{2}{2\alpha+1}} m^{\frac{2\alpha+2}{2\alpha+1}}\right)$ total bits to communicate in order to achieve the centralized optimal rate $O\left(M^{\frac{2}{2\alpha+1}}\left(\frac{\sigma^2}{m}\right)^{\frac{2\alpha}{2\alpha+1}}\right)$. It is tempting to consider the question: Is there a distributed estimator with a total communication budget $O\left(\left(\frac{M}{\sigma}\right)^{\frac{2}{2\alpha+1}} m^{\frac{2\alpha+2}{2\alpha+1}}\right)$ that adaptively achieves the centralized optimal rate over a wide collection of Besov classes $\theta \in \mathcal{B}_{p,q}^{\alpha}(M)$?

To rigorously formulate this problem, let $\tilde{S} \subset (0,\infty) \times (0,\infty) \times (0,\infty] \times (0,\infty]$ be a collection of Besov parameter combinations $(\alpha, M, p, q)$, and $\tilde{C}(\cdot)$ is a function $(0,\infty) \to (0,\infty)$. Let $\mathcal{G}(\tilde{S}, \tilde{C})$ be the set of adaptive distributed estimators that achieve the centralized optimal rate of convergence over Besov classes $\mathcal{B}_{p,q}^{\alpha}(M)$ for all $(\alpha, M, p, q) \in \tilde{S}$. To be precise, $\mathcal{G}(\tilde{S}, \tilde{C})$ is the collection of distributed estimators $\hat{\theta}$ who satisfy the following property: for any $(\alpha, M, p, q) \in \tilde{S}$,

$$\sup_{\theta \in \mathcal{B}_{p,q}^{\alpha}(M)} \mathbb{E}\|\hat{\theta} - \theta\|^2 \leq \tilde{C}(\alpha) M^{\frac{2}{2\alpha+1}} \left(\frac{\sigma^2}{m}\right)^{\frac{2\alpha}{2\alpha+1}}$$

Estimators in $\mathcal{G}(\tilde{S}, \tilde{C})$ are called *statistically-optimal adaptive estimators* over parameter set $\tilde{S}$. We are interested in the minimum expected communication cost among all statistically-optimal adaptive estimators:

$$Q(\tilde{S}, \tilde{C}, \mathcal{B}_{p,q}^{\alpha}(M)) \triangleq \inf_{\hat{\theta} \in \mathcal{G}(\tilde{S}, \tilde{C})} \sup_{\theta \in \mathcal{B}_{p,q}^{\alpha}(M)} \mathbb{E}_\theta L(\hat{\theta})$$

The above quantity, which is called *the minimax communication cost for statistically-optimal adaptive estimators*, serves as a benchmark for the communication-efficiency of estimators in $\mathcal{G}(\tilde{S}, \tilde{C})$. For any statistically-optimal adaptive estimators, its expected communication cost is at least $Q(\tilde{S}, \tilde{C}, \mathcal{B}_{p,q}^{\alpha}(M))$ when estimating a function in $\mathcal{B}_{p,q}^{\alpha}(M)$. The analysis of the minimax communication cost $Q(\tilde{S}, \tilde{C}, \mathcal{B}_{p,q}^{\alpha}(M))$ is divided into two steps: upper bound and lower bound. We first propose in Section 3.1 an adaptive distributed estimator $\hat{\theta}^A$ which can achieve the centralized optimal rate of convergence when $2 \leq p \leq \infty$, and provide a upper bound on the expected communication cost. We then derive in Section 3.2 a lower bound for the rate of convergence of $Q(\tilde{S}_0, \tilde{C}, \mathcal{B}_{p,q}^{\alpha}(M))$ where $\tilde{S}_0$ is collection of all Besov class parameters with $p \geq 2$. The lower bound provides a fundamental limit on the communication cost for a statistically-optimal adaptive estimator, while it matches the upper bound for $\hat{\theta}^A$ on the expected communication cost. Therefore, the proposed distributed estimator $\hat{\theta}^A$ is shown to be the most communication-efficient one among all statistically-optimal adaptive estimators over a wide range of Besov classes.

3.1. *Optimal adaptive procedure by local thresholding.* In order to establish an upper bound on $Q(\tilde{S}, \tilde{C}, \mathcal{B}_{p,q}^{\alpha}(M))$, we first construct a statistically-optimal adaptive distributed procedure which simultaneously achieves the optimal rate of convergence over a wide collection of Besov classes, while the rate of convergence for its expected communication cost matches that of the minimax lower bound given in Section 3.2.

Wavelet thresholding methods have been shown to be a powerful tool for adaptive nonparametric function estimation problems in the conventional centralized settings. Estimators derived from data-driven thresholding rules can automatically adapt to a wide collection of Besov spaces. See Donoho and Johnstone (1995); Abramovich et al. (2006); Cai (1999); Cai and Zhou (2009); Johnstone (2017) and the references therein. However, in the distributed settings, due to the communication constraints, it is typically impossible to estimate individual coordinates accurately by thresholding them all together on the central machine. In such a setting, it is unclear how to optimally threshold on each local machine and efficiently transmit the information to the central machine with minimal communication cost such that a final aggregated estimator is statistically-optimal adaptive. Indeed, it is unclear if this goal is even achievable.

Fortunately, the answer is affirmative. The following "local thresholding" procedure is proposed for adaptive distributed Gaussian sequence estimation. We should emphasize that here "local thresholding" referred to the fact that the thresholding step is carried out on individual local machines, not on the central machine. The meaning is different from that in the conventional wavelet estimation literature in the centralized setting. The general strategy can be summarized as follows. On each local machine, we first select "significant resolution levels" by certain thresholding rule. Only information about the significant resolution levels is transmitted to the central machine, where an estimation subroutine called "ada-MODGAME" is applied to generate good estimates for individual coordinates based on the transcripts collected from the local machines. These estimates will be further processed to yield a final estimate $\hat{\theta}^A$.

Now we are ready to introduce the local thresholding procedure in detail. Let $g : \mathbb{Z} \to \bigcup_{k=1}^{\infty} \{0,1\}^k$ denote the localization encoding function defined in Section 2.1. The estimation procedure is divided into two steps, with the subroutine *ada-MODGAME* in the second step of the procedure.

**First step: Generate the transcripts on the local machines by thresholding.** For $1 \leq i \leq m$, on the $i$-th machine:

1. Define the set of "significant resolution levels" on the $i$-th machine by

$$J_i = \{0, 1, 2, ..., (\lfloor 2\log m \rfloor)\} \bigcup \{j \geq \lfloor 2\log m \rfloor + 1 : \sum_{k=1}^{n_j} X_{i,jk}^2 \geq n_j \sigma^2 (1 + \frac{\Lambda_1}{m})\},$$

where $\Lambda_1 > 0$ is a prespecified parameter. Only those coordinates at the resolution levels in the set $J_i$ are processed as part of the transcript outputs from the $i$-th machine. All the resolution levels that are not in $J_i$ are considered to be "locally thresholded", because the signal strength on those resolution levels is weak.

2. If $i = 1$, the output transcript $Z_1$ is the collection of the "crude localization" strings $Z_{1,jk}$, $j \in J_1, k \in [n_j]$ where $Z_{1,jk}$ is defined as

$$Z_{1,jk} = g(\lfloor X_{1,jk}/\sigma \rfloor);$$

If $2 \leq i \leq 1 + \lfloor \log^2 m \rfloor$, the output transcript $Z_i$ is the collection of the "finer localization" strings $Z_{i,jk}$, $j \in J_i, k \in [n_j]$ where $Z_{i,jk}$ is defined as

$$Z_{i,jk} = g(\lfloor X_{i,jk}/\sigma \rfloor \bmod \lfloor \log m \rfloor);$$

If $i \geq 2 + \lfloor \log^2 m \rfloor$ the output transcript $Z_i$ is the collection of the "refinement" strings $Z_{i,jk}$, $j \in J_i, k \in [n_j]$ where $Z_{i,jk}$ is defined as

$$Z_{i,jk} = \lfloor X_{i,jk}/\sigma \rfloor \bmod 8.$$

**Second step: Generate the distributed estimator $\hat{\theta}$ on the central machine.** The central machine receives the transcripts $Z_1, Z_2, ..., Z_m$ from the local machines. Note that the code words in $Z_1, Z_2, ..., Z_m$ are all uniquely decodable, thus the central machine is able to recover short strings $Z_{i,jk}$ for $i \in [m], j \in J_i, k \in [n_j]$. Also, note that the total number of short strings from the $i$-th machine is $\sum_{j \in J_i} 2^j$, so from the binary representation of the total number of short strings from the $i$-th machine, one can recover significant resolution level $J_i$.

To warp up, from those transcripts that the central machine receives

- significant resolution levels on the local machines $J_1, J_2, ..., J_m$.
- short strings $Z_{i,jk}$ for $i \in [m], j \in J_i, k \in [n_j]$.

Let $\hat{J}$ be defined as

$$\hat{J} \triangleq \{j : j \in J_1; \sum_{i=2}^{1+\lfloor \log^2 m \rfloor} \mathbb{I}_{\{j \in J_i\}} \geq \frac{\lfloor \log^2 m \rfloor}{2}; \sum_{i=2+\lfloor \log^2 m \rfloor}^{m} \mathbb{I}_{\{j \in J_i\}} \geq \frac{m - 1 - \lfloor \log^2 m \rfloor}{2}\}$$

Intuitively, $\hat{J}$ is the set of resolution levels that are significant on most local machines. The resolution levels within $\hat{J}$ will be estimated whereas those not in $\hat{J}$ will be zero out (thresholded).

The final estimator $\hat{\theta}^A$ is constructed as follows: For $j = 1, 2, ...$,

- If $j \notin \hat{J}$, let
$$\hat{\theta}^A_{jk} = 0 \text{ for all } k \in [n_j].$$

- If $j \leq \lfloor 2 \log m \rfloor$, let $S_j = [m]$ and

$$(\hat{\theta}^*_{j1}, \hat{\theta}^*_{j2}, ..., \hat{\theta}^*_{jn_j}) = \hat{f}_{\text{ada}}(S_j, \{Z_{i,jk} : i \in S_j, k \in [n_j]\})$$

  be the output of the subroutine "ada-MODGAME". Then apply the thresholding rule to get the final estimate

$$(\hat{\theta}^A_{j1}, \hat{\theta}^A_{j2}, ..., \hat{\theta}^A_{jn_j}) = \begin{cases} (\hat{\theta}^*_{j1}, \hat{\theta}^*_{j2}, ..., \hat{\theta}^*_{jn_j}) & \text{if } \sum_{k=1}^{n_j} (\hat{\theta}^*_{jk})^2 \geq \Lambda_2 \frac{n_j \sigma^2}{m} \\ (0, 0, 0, ..., 0) & \text{otherwise} \end{cases}$$

  where $\Lambda_1 > 0$ is a prespecified parameter.

- If $j \geq \lfloor 2 \log m \rfloor + 1$ and $j \in \hat{J}$, define $S_j = \{i \in [m] : j \in J_i\}$, and let

$$(\hat{\theta}^A_{j1}, \hat{\theta}^A_{j2}, ..., \hat{\theta}^A_{jn_j}) = \hat{f}_{\text{ada}}(S_j, \{Z_{i,jk} : i \in S_j, k \in [n_j]\})$$

  be the output of the subroutine "ada-MODGAME".

**Subroutine: ada-MODGAME**

Input: $\sigma, m, j, n_j, S_j, \{Z_{i,jk} : i \in S_j, k \in [n_j]\}$.

For each $k \in [n_j]$, do following steps:

1. Because $g(x)$ is uniquely decodable, from $Z_{1,jk} = g(\lfloor X_{i,jk}/\sigma \rfloor)$ one can recover the value of $\lfloor X_{i,jk}/\sigma \rfloor$. Let $I^a_{jk}$ be a left-closed-right-open interval of length $m$ defined as

$$I^a_{jk} \triangleq \begin{cases} \left[\lfloor X_{i,jk}/\sigma \rfloor - \frac{\lfloor \log m \rfloor - 1}{2}, \lfloor X_{i,jk}/\sigma \rfloor + \frac{\lfloor \log m \rfloor + 1}{2}\right) & \text{if } \lfloor \log m \rfloor \text{ is an odd number} \\ \left[\lfloor X_{i,jk}/\sigma \rfloor - \frac{\lfloor \log m \rfloor}{2}, \lfloor X_{i,jk}/\sigma \rfloor + \frac{\lfloor \log m \rfloor}{2}\right) & \text{if } \lfloor \log m \rfloor \text{ is an even number} \end{cases}.$$

2. Let $S^b_j \triangleq S_j \cap \{i : 2 \leq i \leq \lfloor \log^2 m \rfloor + 1\}$ be the set of machines that output the finer localization strings. Let $z^b_{ik} \triangleq \text{argmax}_{z'} \sum_{i \in S^b_j} \mathbb{I}_{\{Z_{i,jk} = z'\}}$ be the mode statistic among $Z_{i,jk}, i \in S^b_j$. Note that the length of $I^a_{jk}$ is $\lfloor \log m \rfloor$, so there will be exactly one integer $x^b_{jk} \in I^a_{jk}$ satisfying

$$x^b_{jk} \bmod \lfloor \log m \rfloor = g^{-1}(z^b_{ik}).$$

  Let $I^b_{jk}$ be an interval of length 3 defined by

$$I^b_{jk} \triangleq [x^b_{jk} - 1, x^b_{jk} + 1].$$

3. Let $S^h_j \triangleq S_j \cap \{i : i \geq \lfloor \log^2 m \rfloor + 2\}$ be the set of machines that output the refinement strings. Let $p^h$ be the proportion of those refinement strings whose value is equal to $g(x^b_{jk} - 2 \bmod 8)$:

$$p^h \triangleq \text{card}(S^h_j)^{-1} \sum_{i \in S^h_j} (\mathbb{I}_{\{Z_{i,jk} = g(x^b_{jk} - 2 \bmod 8)\}}$$

Define a function
$$h_{jk}(y) \triangleq \sum_{l=-\infty}^{\infty} \int_{x^b_{jk}-2+8l}^{x^b_{jk}-1+8l} \phi_1(x-y)dx.$$

It is easy to see that $h_{jk}(y)$ is a strictly decreasing function on $I^b_{jk}$. Let $h^{-1}_{jk}(y)$ be the inverse function of $h_{jk}(y)$ which maps $h_{jk}(I^b_{jk})$ to $I^b_{jk}$. Finally the estimate can be calculated by

$$\hat{\theta}^*_{jk} = \begin{cases} (x^b_{jk}+1)\sigma & \text{if } p^h \leq h_{jk}(x^b_{jk}+1) \\ h^{-1}_{jk}(p^h)\sigma & \text{if } h_{jk}(x^b_{jk}+1) < p^h < h_{jk}(x^b_{jk}-1) \\ (x^b_{jk}-1)\sigma & \text{if } p^h \geq h_{jk}(x^b_{jk}-1) \end{cases}.$$

Output: $\hat{\theta}^*_{jk}$ for $k \in [n_j]$.

We have given above a detailed construction of the local thresholding estimator $\hat{\theta}^A$. The following theorem provides a theoretical guarantee for the statistical performance and communication cost of the proposed procedure over the Besov classes $\mathcal{B}^\alpha_{p,q}(M)$ with $\alpha > 0, M \geq \sigma, 1 < q \leq \infty$, and $2 \leq p \leq \infty$.

THEOREM 3 (Upper Bound for the Communication Cost). *If $\Lambda_1 > 10$ and $\Lambda_2$ is chosen sufficiently large, there exists a constant $C > 0$ such that, the local thresholding estimator $\hat{\theta}^A$ is adaptively rate-optimal, i.e.*

$$\sup_{\theta \in \mathcal{B}^\alpha_{p,q}(M)} \mathbb{E}\|\hat{\theta}^A - \theta\|^2 \leq C M^{\frac{2}{2\alpha+1}} \left(\frac{\sigma^2}{m}\right)^{\frac{2\alpha}{2\alpha+1}}$$

*and we also have*

$$\sup_{\theta \in \mathcal{B}^\alpha_{p,q}(M)} \mathbb{E}_\theta L(\hat{\theta}^A) \leq C \left(m^3 + \left(\frac{M}{\sigma}\right)^{\frac{2}{2\alpha+1}} m^{\frac{2\alpha+2}{2\alpha+1}}\right)$$

*for all $\alpha > 0, M \geq \sigma, 1 < q \leq \infty$, and $2 \leq p \leq \infty$.*

REMARK 4. The proof of Theorem 3 is involved due to the fact that, after thresholding on the local machines, the conditional distribution of the observations given that their resolution level is selected into the significant set $J_i$ is no longer Gaussian. Lemma **??** (from the supplementary material Cai and Wei (2020b)) is the key to the proof, which shows that the ada-MODGAME subroutine is robust even if the additive noise is slightly different from Gaussian distribution.

REMARK 5. One of the merits of the local thresholding estimator $\hat{\theta}^A$ is its "communication-adaptivity", which means the communication cost of the estimation procedure is also adaptive to the smoothness of the underlying function. Compared to the two-point adaptive procedure proposed in the previous work Szabo and van Zanten (2020) which is able to achieve adaptation with smoothness less than certain threshold, our newly proposed local thresholding procedure requires no prior knowledge on the range of the smoothness parameters, and is able to achieve statistical adaptation over a wide collection of Besov classes. The user can apply local thresholding procedure to obtain adaptation over the Besov classes $\mathcal{B}_{p,q}^{\alpha}(M)$ as long as $p \geq 2$ with guaranteed minimum communication cost.

3.2. *Lower bound analysis.* In this subsection, we are going to obtain a lower bound for the minimax communication cost for statistically-optimal adaptive estimators, which is instrumental in establishing the optimal rate of convergence. Before we establish a lower bound for the minimax communication cost $Q(\tilde{S}, \tilde{C}, \mathcal{B}_{p,q}^{\alpha}(M))$, we first state the following theorem, which gives a lower bound for the communication cost when the estimator achieves statistical-optimal rate of convergence in two different Besov classes.

THEOREM 4 (Lower bound for communication cost for two-point adaptation). *For any distributed estimator $\hat{\theta}$, let $\mathcal{B}_{p_1,q_1}^{\alpha_1}(M_1)$ and $\mathcal{B}_{p_2,q_2}^{\alpha_2}(M_2)$ be two different Besov classes. If there exists a constant $C > 0$ such that $M_1 \leq C\sigma m^{2\alpha_1 + \frac{1}{2}}$, and*

$$(9) \qquad \sup_{\theta \in \mathcal{B}_{p_l,q_l}^{\alpha_l}(M_l)} \mathbb{E}\|\hat{\theta} - \theta\|^2 \leq C M_l^{\frac{2}{2\alpha_l+1}} \left(\frac{\sigma^2}{m}\right)^{\frac{2\alpha_l}{2\alpha_l+1}} \text{ for } l = 1, 2.$$

*Then there exists a constant $c > 0$ (depending on $C$) such that*

$$\sup_{\theta \in \mathcal{B}_{p_2,q_2}^{\alpha_2}(M_2)} \mathbb{E}L(\hat{\theta}) \geq c \left( \left(\frac{M_1}{\sigma}\right)^{\frac{2}{2\alpha_1+1}} m^{\frac{2\alpha_1+2}{2\alpha_1+1}} + \left(\frac{M_2}{\sigma}\right)^{\frac{2}{2\alpha_2+1}} m^{\frac{2\alpha_2+2}{2\alpha_2+1}} \right).$$

REMARK 6. If one sets $\sigma = \sqrt{m/n}$, $M_1 = M_2 = 1$ and $\alpha_2 > \alpha_1 > \frac{\log n}{4 \log m} - \frac{1}{2}$, the above Theorem 4 recovers the result of Theorem 2.4 in Szabo and van Zanten (2020) which shows that two-point adaptation is impossible without additional communication cost when $m^{4\alpha+2} \gg n$. Comparing with the previous result, the result given in Theorem 4 here is stronger because we prove the lower bound for the communication cost $\sup_{\theta \in \mathcal{B}_{p_2,q_2}^{\alpha_2}(M_2)} \mathbb{E}L(\hat{\theta})$ under the only assumption that $\hat{\theta}$ is adaptive. In particular, no upper bound is imposed on $\sup_{\theta \in \mathcal{B}_{p_1,q_1}^{\alpha_1}(M_1)} \mathbb{E}_\theta L(\hat{\theta})$, which is in fact necessary to obtain Theorem 2.4 in Szabo and van Zanten (2020).

The above Theorem 4 only considers two-point adaptation between two specific Besov classes. However, in real data application, we are more interested in developing estimators that are able to adapt to a wide range of parameter spaces, such as our adaptive estimator $\hat{\theta}^A$. It is necessary to extend the above Theorem 4 to a general lower bound on $Q(\tilde{S}, \tilde{C}, \mathcal{B}_{p,q}^{\alpha}(M))$.

We define $\tilde{S}_0 = \{(\alpha, M, p, q) : \alpha > 0, M \geq \sigma, 2 \leq p \leq \infty, 1 < q \leq \infty\}$ a wide collection of Besov class parameters. The following lower bound on $Q(\tilde{S}_0, \tilde{C}, \mathcal{B}_{p,q}^{\alpha}(M))$ shows a fundamental limit on the communication cost of statistically-optimal estimators over Besov classes $\mathcal{B}_{p,q}^{\alpha}(M)$ where $(\alpha, p, q, M) \in \tilde{S}_0$. In view of the upper bound to be given in Section 3.1 that is achieved by the adaptive distributed estimator $\hat{\theta}^A$, the lower bound is rate optimal.

THEOREM 5 (Lower bound for the communication cost over Besov ball collection $\tilde{S}_0$). *For any $\tilde{C} : (0, \infty) \to (0, \infty)$ and $(\alpha, M, p, q) \in \tilde{S}_0$, there exists a constant $c > 0$ such that*

$$(10) \qquad Q(\tilde{S}_0, \tilde{C}, \mathcal{B}_{p,q}^{\alpha}(M)) \geq c \left( m^3 + \left(\frac{M}{\sigma}\right)^{\frac{2}{2\alpha+1}} m^{\frac{2\alpha+2}{2\alpha+1}} \right).$$

REMARK 7. The lower bound in Theorem 5 shows that, if a distributed estimator adaptively achieves the optimal rate of convergence over the all Besov classes where $p \geq 2$, the minimum required expected communication cost for estimating functions in $\mathcal{B}_{p,q}^{\alpha}(M)$ is of order $m^3 + \left(\frac{M}{\sigma}\right)^{\frac{2}{2\alpha+1}} m^{\frac{2\alpha+2}{2\alpha+1}}$. The additional communication cost, which is of order $m^3$ and not depending on the values of $\alpha, M, p, q$ and $\sigma$, is required and necessary for constructing an adaptive estimator. When $m \gtrsim \left(\frac{M}{\sigma}\right)^{\frac{2}{4\alpha+1}}$, the cost of adaptation is significant.

REMARK 8. Although in Theorem 5 we provide a lower bound on $Q(\tilde{S}, \tilde{C}, \mathcal{B}_{p,q}^{\alpha}(M))$ where $\tilde{S} = \tilde{S}_0$, the same lower bound also holds when $\tilde{S}$ is other sufficiently large Besov ball collections. With the help from Theorem 4, we are able to establish lower bounds for other Besov ball collection $\tilde{S}$.

REMARK 9. The techniques used to prove Theorems 4 and 5 can be of independent interest. Roughly speaking, if the algorithm aims to perform well on both $B_{p_1,q_1}^{\alpha_1}(M_1)$ and $B_{p_2,q_2}^{\alpha_2}(M_2)$ where $\alpha_1 < \alpha_2$, since we cannot tell whether each local sample is drawn from $B_{p_1,q_1}^{\alpha_1}(M_1)$ or $B_{p_2,q_2}^{\alpha_2}(M_2)$ on the local machines, the algorithm needs to transmit more bits than non-adaptive estimation for $B_{p_2,q_2}^{\alpha_2}(M_2)$, because it also needs to estimate well in $B_{p_1,q_1}^{\alpha_1}(M_1)$. More specifically, we prove that the local machines cannot "distinguish" samples that is drawn from a null model ($\theta = \vec{0}$) or drawn from a mixture of models with $\theta$ having $m^2$ non-zero elements. If the observations are truly drawn from the mixture, the minimum communication cost required to achieve the statistical optimal rate of convergence is of order $m^3$. Thus one can further show that the minimax communication cost is at least $\Omega(m^3)$ even if $\theta = \vec{0}$. This is a key step in the argument for establishing Theorems 4 and 5.

A similar technique was also used in Szabo and van Zanten (2020). But a finer analysis is needed here, especially for the key Claim 3 where we first prove a *conditional strong* *data processing inequality* and use it to establish a stronger result without unnecessary assumptions.

LEMMA 1 (Conditional strong data processing inequality). *For $t > 0$ and $k \in \mathbb{Z}_+$, let $\theta$ be a random vector uniformly distributed on the set $\{-t\sigma, t\sigma\}^k$ and let $X \sim N(\theta, \sigma^2 I_k)$. Let $D \subseteq \mathbb{R}^k$ be a $k$-dimensional region such that the event $X \in D$ is independent with $\theta$ and let $Z$ be a random variable such that $\theta \to X \to Z$ forms a Markov chain. Then*

$$I(Z;\theta|X \in D)\mathbb{P}(X \in D) \leq 256t^2(H(Z|X \in D)\mathbb{P}(X \in D) + H(\{X \in D\})),$$

*where $I(\cdot;\cdot|\cdot)$, $H(\cdot)$, and $H(\cdot|\cdot)$ denote conditional mutual information, entropy, and conditional entropy respectively.*

The definitions of the conditional mutual information $I(\cdot;\cdot|\cdot)$, entropy $H(\cdot)$, and conditional entropy $H(\cdot|\cdot)$ are given in Section 6.1. Note that the classical strong data processing inequality for the Gaussian channels serves as a special case if we set $D = \mathbb{R}^k$. The above inequality is the key to the proof of Theorem 4. We omit the proof of Lemma 1 since it is similar to the proof of Claim 3 in the proof of Theorem 4.

The upper and lower bounds given in Theorems 3 and 5 together establish the minimax rate of communication cost for statistically-optimal adaptive estimators:

$$(11) \qquad Q(\tilde{S}_0, \tilde{C}, \mathcal{B}_{p,q}^\alpha(M)) \asymp m^3 + \left(\frac{M}{\sigma}\right)^{\frac{2}{2\alpha+1}} m^{\frac{2\alpha+2}{2\alpha+1}}$$

where $\tilde{C}$ is large enough and recall that $\tilde{S}_0 = \{(\alpha, M, p, q) : \alpha > 0, M > \sigma, 2 \leq p \leq \infty, 1 < q \leq \infty\}$. The minimax rate (11) also implies that $\hat{\theta}^A$ is the optimal adaptive distributed estimator with respect to both statistical performance and communication cost.

1. The estimator $\hat{\theta}^A$ simultaneously achieves the centralized optimal rate over the Besov classes $\mathcal{B}_{p,q}^\alpha(M)$ for all $\alpha > 0, M \geq \sigma$, $1 < q \leq \infty$, and $2 \leq p \leq \infty$. There is no statistical cost of adaptation in terms of the rate of convergence.
2. Among all the statistically-optimal adaptive estimators, the expected communication cost for $\hat{\theta}^A$ is rate-optimal over the Besov classes $\mathcal{B}_{p,q}^\alpha(M)$ for all $\alpha > 0, M \geq \sigma$, $1 < q \leq \infty$, and $2 \leq p \leq \infty$.

REMARK 10. Compared with the minimum communication cost $\left(\frac{M}{\sigma}\right)^{\frac{2}{2\alpha+1}} m^{\frac{2\alpha+2}{2\alpha+1}}$ for achieving the optimal rate of convergence in the minimax setting in (8), an additional communication cost of order $m^3$ bits is needed to achieve the adaptation over a collection of Besov classes. The term $m^3$ can be viewed as the communication cost of adaptation. This interplay between communication and statistical adaptation in the distributed setting is an interesting phenomenon: It costs more bits to communicate in order to achieve adaptivity. In contrast, statistical adaptation can be achieved for free in the centralized setting (Donoho and Johnstone, 1995; Johnstone, 2017).

**4. Numerical Studies.** The proposed seq-MODGAME estimator $\hat{\theta}^O$ and the adaptive local thresholding estimator $\hat{\theta}^A$ are easily to implement. In this section, we conduct simulation studies to investigate the numerical performance of these two estimators in various settings.

4.1. *The seq-MODGAME estimator $\hat{\theta}^O$.* We first study the seq-MODGAME estimator $\hat{\theta}^O$ proposed in Section 2. We generate i.i.d data according to the distributed Gaussian sequence model (3) on $m = 100$ different virtual machines, where the mean vector $\theta$ is the wavelet coefficients of certain specified underlying function. The underlying function $f$ is chosen as

$$f(t) = \sin(4\pi t) + 0.7\cos(18\pi t) \quad t \in [0,1]$$

and the noise level $\sigma = 1/16$.

We apply the optimal seq-MODGAME estimator $\theta^O$ to estimate wavelet coefficients of $f$ given their noisy observations on virtual machines. Afterwards, we transform estimated wavelet coefficients back to estimated smooth functions $\hat{f}^O$. The results are shown in Figure 1. As more and more bits are allowed to communicate, the mean squared error are decreasing so that the estimate is becoming more and more accurate.

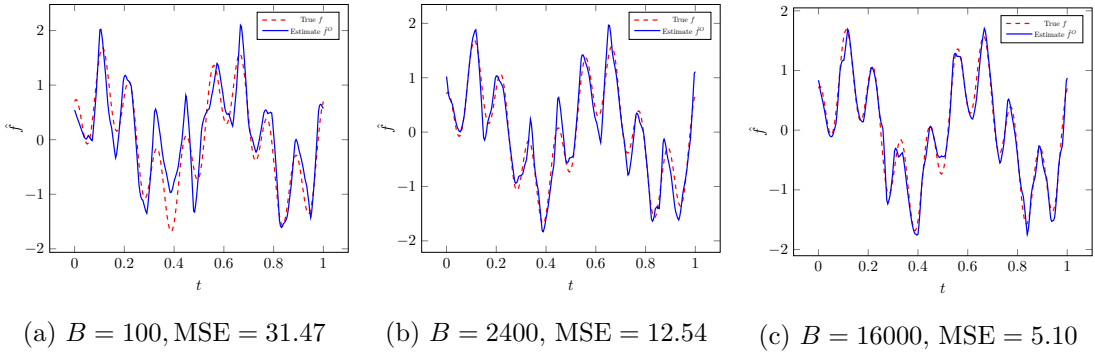

(a) $B = 100, \text{MSE} = 31.47$　　(b) $B = 2400, \text{MSE} = 12.54$　　(c) $B = 16000, \text{MSE} = 5.10$

Fig 1: Estimate given by the optimal seq-MODGAME estimator $\hat{\theta}^O$ under the communication constraints. For different choices of total communication budgets $B = 100, 2400, 16000$, we illustrate an example of estimated function $\hat{f}^O$ in each figure. The mean squared error through 1000 trials are also given below each figure.

4.2. *The local thresholding estimator $\hat{\theta}^A$.* Similar to the setting in Section 4.1, we generate i.i.d data according to the distributed Gaussian sequence model (3) and set $m = 100, \sigma = 1/16$. However, in this simulation study we work on three different choices for the underlying functions $f = f_1, f_2$ or $f_3$:

$$f_1(t) = 1.5\sin(4\pi t) \qquad t \in [0, 1];$$
$$f_2(t) = \sin(4\pi t) + 0.7\cos(18\pi t) \qquad t \in [0, 1];$$
$$f_3(t) = 0.8\sin(4\pi t) + 0.5\cos(18\pi t) + 0.5\cos(44\pi t) \qquad t \in [0, 1].$$

The three functions given above are designed to have different smoothness. $f_1$ is the smoothest function among the three functions whereas $f_3$ is the most "wiggly" one. We expect to see a data-driven estimator can adapt to their smoothness automatically during the estimation.

Similarly, given random distributed data generated by adding noise to the wavelet coefficients of $f_1, f_2$ and $f_3$ respectively, we apply the local thresholding estimator $\hat{\theta}^A$ to estimate the wavelet coefficients. The estimated smooth functions $\hat{f}^A$ are obtained by reversed discrete wavelet transform on the estimated wavelet coefficients. The results are shown in Figure 2. It can be clearly seen from simulation that, when the underlying function are relatively smooth, the local thresholding estimator requires less communication cost while achieves better statistical accuracy. The numerical results are consistent with the theory, which shows the local thresholding estimator can adapt to the smoothness of the underlying function.

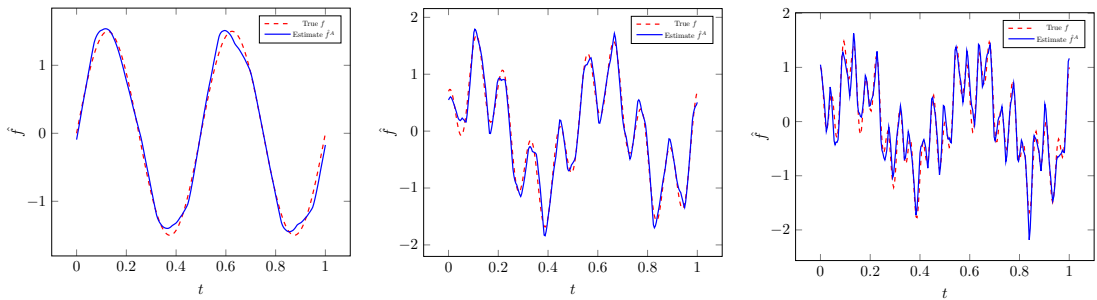

(a) $f_1 : \mathbb{E}L = 3330, \text{MSE} = 2.06$ (b) $f_2 : \mathbb{E}L = 8083, \text{MSE} = 5.03$ (c) $f_3 : \mathbb{E}L = 15862, \text{MSE} = 8.9$

Fig 2: Estimate given by the local thresholding estimator $\hat{\theta}^A$. Under different choices of ground truth functions $f_1, f_2, f_3$, we illustrate an example of estimated function $\hat{f}^A$ in each figure. The expected communication cost and their mean squared error through 1000 trials are also given below each figure.

**5. Discussion.** In the present paper, both distributed minimax and distributed adaptive estimation under the communication constraints were studied for the Gaussian sequence model and white noise model. Optimal minimax rate of convergence is established and the cost of adaptation is characterized. In addition, a data-driven adaptive distributed estimator with theoretical guarantees is constructed. Several technical tools and the for-mulation for the study of the interplay between adaptation and communication cost can be of independent interest.

Distributed nonparametric function estimation is still very much a new area with a range of interesting open problems. One such problem is the construction of an adaptive distributed procedure for Gaussian sequence estimation under a fixed communication constraint. It is notable that the communication cost for the local thresholding procedure $\hat{\theta}^A$ is related to the smoothness of the underlying function. When the communication budget is tight, there is not enough budget to implement the local thresholding procedure. Therefore, it will be useful to have an estimator whose communication cost is controlled, while its estimation accuracy is adaptive to the smoothness of underlying function.

In the present paper, we focused on estimation over the Besov classes with $p \geq 2$. Another direction is the study of distributed Gaussian sequence estimation over the Besov classes with $p < 2$. Similar to the centralized setting, the case $p < 2$ is very different from the case $p \geq 2$ in the distributed setting. The techniques developed in the present paper are not sufficient for the case $p < 2$ and we leave this case for future work.

Besides the white noise model considered in the present paper, it is also interesting to study other related nonparametric function estimation problems, including nonparametric density estimation, nonparametric regression with fixed design, and nonparametric regression with random design, which have all been well studied in the centralized setting. In particular, it is shown that these three models are asymptotic equivalent to the white noise model (Nussbaum, 1996; Brown and Low, 1996; Brown et al., 2002, 2004) in the centralized setting under mild regularity conditions when the smoothness parameter $\alpha > \frac{1}{2}$. Practically, for example, by applying the root-unroot algorithm to the binned data (Brown et al., 2010), the density estimation problem can essentially be turned into the problem of nonparametric regression with fixed design. However, in the distributed settings, these four problems may exhibit different asymptotic behaviors due to the communication constraints. In the distributed setting, nonparametric density estimation, nonparametric regression with fixed design, and nonparametric regression with random design merit careful and separate investigations. We leave them for future work.

Broadly speaking, virtually any problem studied in the classical centralized setting has its counterpart in the distributed setting. Examples include minimax and adaptive estimation of linear and quadratic functionals as well as hypothesis testing under these nonparametric function models. It is challenging to develop a general optimality theory and construct statistically optimal distributed procedures under the communication constraints, New technical tools for both the lower bound and upper bound analyses are needed.

**6. Proofs.** We prove Theorems 4 and 5 in this section. For reasons of space, the proofs of the other theorems, propositions and additional technical lemmas are given in the supplementary material (Cai and Wei, 2020b).

6.1. *Notation and definitions.* For any finite $S$, denote $\mathbb{U}(A)$ be a uniform distribution on $S$. For any $a, b$, let $a \lesssim b$ denote there exists a universal constant $C > 0$ such that $a \leq Cb$, whereas $a \gtrsim b$ denotes there exists a universal constant $c > 0$ such that $a \geq cb$. For any discrete random variables $X, Y$ supported on $\mathcal{X}, \mathcal{Y}$, the entropy $H(X)$, conditional entropy $H(X|Y)$, and mutual information $I(X; Y)$ are defined as

$$H(X) \triangleq -\sum_{x \in \mathcal{X}} \mathbb{P}(X = x) \log \mathbb{P}(X = x),$$

$$H(X|Y) \triangleq -\sum_{x \in \mathcal{X}, y \in \mathcal{Y}} \mathbb{P}(X = x, Y = y) \log \mathbb{P}(X = x|Y = y),$$

$$I(X; Y) \triangleq \sum_{x \in \mathcal{X}, y \in \mathcal{Y}} \mathbb{P}(X = x, Y = y) \log \frac{\mathbb{P}(X = x|Y = y)}{\mathbb{P}(X = x)}.$$

6.2. *Proof of Theorem 4.* It follows from Theorem 2 that for any estimator $\hat{\theta}$ satisfying $\sup_{\theta \in \mathcal{B}^\alpha_{p,q}(M)} \mathbb{E}_\theta L(\hat{\theta}) \leq B$, we have

$$\sup_{\theta \in \mathcal{B}^\alpha_{p,q}(M)} \mathbb{E}\|\hat{\theta} - \theta\|^2 \geq cM^{\frac{2}{\alpha+1}} \left(\frac{\sigma^2}{B}\right)^{\frac{\alpha}{\alpha+1}}$$

for some constant $c > 0$. By the assumption,

$$\sup_{\theta \in \mathcal{B}^{\alpha_2}_{p_2,q_2}(M_2)} \mathbb{E}\|\hat{\theta} - \theta\|^2 \leq CM_2^{\frac{2}{2\alpha_2+1}} \left(\frac{\sigma^2}{m}\right)^{\frac{2\alpha_2}{2\alpha_2+1}}.$$

So it follows that

$$\sup_{\theta \in \mathcal{B}^{\alpha_2}_{p_2,q_2}(M_2)} \mathbb{E}_\theta L(\hat{\theta}) \gtrsim \left(\frac{M_2}{\sigma}\right)^{\frac{2}{2\alpha_2+1}} m^{\frac{2\alpha_2+2}{2\alpha_2+1}}.$$

To prove Theorem 4, it now suffices to show

$$\sup_{\theta \in \mathcal{B}^{\alpha_2}_{p_2,q_2}(M_2)} \mathbb{E}_\theta L(\hat{\theta}) \gtrsim \left(\frac{M_1}{\sigma}\right)^{\frac{2}{2\alpha_1+1}} m^{\frac{2\alpha_1+2}{2\alpha_1+1}}.$$

The remain part of the proof aims to prove the above inequality.

Define the constant $\lambda$ (only depends on $C$) and variable $u$ as follows:

$$\lambda = \max\{10, 32\sqrt{C}\},$$

$$u = \left(\frac{M_1}{\sigma}\right)^{\frac{2}{2\alpha_1+1}} m^{\frac{1}{2\alpha_1+1}}.$$

Define the set of sequences

$$S_{m,u} \triangleq \left\{ \left(\tau_1 \frac{\lambda\sigma}{\sqrt{m}}, \tau_2 \frac{\lambda\sigma}{\sqrt{m}}, ..., \tau_u \frac{\lambda\sigma}{\sqrt{m}}, 0, 0, ...\right) : \tau_1, \tau_2, ..., \tau_u \in \{-1, +1\} \right\}.$$

Since for any $\theta \in S_{m,u}$ and $p_1, q_1 < \infty$ we have

$$|\theta|_{b_{p_1,q_1}^{\alpha_1}} = \left( \sum_{j=0}^{\infty} \left( 2^{j(\alpha_1+1/2-1/p_1)} \left( \sum_{k=1}^{2^j} |\theta_{jk}|^{p_1} \right)^{1/p_1} \right)^{q_1} \right)^{1/q_1}$$

$$\leq \left( \sum_{j=0}^{\lfloor \log u \rfloor + 1} 2^{jq_1(\alpha_1+1/2)} \right)^{1/q} \frac{\lambda\sigma}{\sqrt{m}}$$

$$\leq \left( \frac{(2u)^{q_1(\alpha_1+1/2)}}{1 - 2^{-q_1(\alpha_1+1/2)}} \right)^{1/q_1} \frac{\lambda\sigma}{\sqrt{m}}$$

$$= \lambda \left( \frac{2^{q_1(\alpha_1+1/2)}}{1 - 2^{-q_1(\alpha_1+1/2)}} \right)^{1/q_1} u^{\alpha_1+1/2} \frac{\sigma}{\sqrt{m}} \leq M_1.$$

When $p_1 = \infty$ or $q_1 = \infty$, the above inequality also holds by similar argument. Therefore we have $S_{m,u} \subset \mathcal{B}_{p_1,q_1}^{\alpha_1}(M_1)$.

Since we have assumed

$$\sup_{\theta \in \mathcal{B}_{p_1,q_1}^{\alpha_1}(M_1)} \mathbb{E}\|\hat{\theta} - \theta\|^2 \leq C M_1^{\frac{2}{2\alpha_1+1}} \left( \frac{\sigma^2}{m} \right)^{\frac{2\alpha_1}{2\alpha_1+1}}.$$

Note that the maximum risk is lower bounded by the Bayesian risk, assign to $\theta$ a uniform prior $\theta \sim \mathbb{U}(S_{m,u})$, then we have

$$\mathbb{E}_{\theta \sim \mathbb{U}(S_m)}\|\hat{\theta} - \theta\|^2 \leq CM_1^{\frac{2}{2\alpha_1+1}} \left( \frac{\sigma^2}{m} \right)^{\frac{2\alpha_1}{2\alpha_1+1}} = Cu\frac{\sigma^2}{m}.$$

In the following proof, we are going to provide several claims and prove each claim accordingly. Let $Q_0$ denote the probability law of $X_1$ when $\theta = (0, 0, 0, ..., 0, ...)$. Let $Q_m$ denote the probability law of $X_1$ when $\theta \sim \mathbb{U}(S_{m,u})$. Note that there are multiple distributions we need to consider, we shorthand the probability, expectation, entropy and mutual information when $\theta = (0, 0, 0, ..., 0, ...)$ as $\mathbb{P}_0, \mathbb{E}_0, H_0$ and $I_0$ respectively. Similarly we use shorthands $\mathbb{P}_m, \mathbb{E}_m, H_m$ and $I_m$ to denote those quantities when $\theta \sim \mathbb{U}(S_{m,u})$

CLAIM 1. *We have $I_m(\hat{\theta}, \theta) \geq \frac{15}{16}u$.*

PROOF OF CLAIM 1: Define $\hat{\theta}^* \triangleq \mathcal{P}_{S_{m,u}}(\hat{\theta})$ be the nearest point in $S_{m,u}$ to $\hat{\theta}$. Then we have

$$\mathbb{E}_m \|\hat{\theta}^* - \theta\|^2 \leq 4 \mathbb{E}_m \|\hat{\theta} - \theta\|^2 \leq 4Cu\frac{\sigma^2}{m}. \tag{12}$$

Note that $\hat{\theta}^* \in S_m$ thus we can reparametrize $\hat{\theta}^*$ to

$$\hat{\theta}^* = (\hat{\tau}_1 \frac{\lambda\sigma}{\sqrt{m}}, \hat{\tau}_2 \frac{\lambda\sigma}{\sqrt{m}}, ..., \hat{\tau}_{m^2} \frac{\lambda\sigma}{\sqrt{m}}, 0, 0, ...) \quad \text{where } \hat{\tau}_1, \hat{\tau}_2, ..., \hat{\tau}_{m^2} \in \{-1, +1\}$$

Then we can simplify (12) to

$$\mathbb{E}_m \sum_{k=1}^{m^2} (\hat{\tau}_k - \tau_k)^2 \leq 4C\lambda^{-2} u. \tag{13}$$

Recall that $\lambda = \max\{10, 32\sqrt{C}\}$. Substitute into (13) we have

$$\mathbb{E}_m \sum_{k=1}^{m^2} (\hat{\tau}_k - \tau_k)^2 \leq \frac{1}{256} u.$$

Apply Fano's inequality, we can conclude

$$\sum_{k=1}^{m^2} H_m(\tau_k | \hat{\tau}_k) \leq \frac{1}{16} u.$$

The following lemma is instrumental to establish later results:

LEMMA 2. *If $A$ is a random variable and $Y_1, Y_2, ..., Y_d$ are independent random variables, then*

$$I(A; (Y_1, Y_2, ..., Y_d)) \geq \sum_{k=1}^{d} I(A; Y_k).$$

Note that $\tau_1, \tau_2, ..., \tau_k$ are i.i.d Rademacher variables, apply Lemma 2 we have

$$I_m(\hat{\theta}^*; \theta) = I_m(\hat{\theta}^*; (\tau_1, \tau_2, ..., \tau_{m^2})) \geq \sum_{k=1}^{m^2} I_m(\hat{\theta}^*; \tau_k) \geq \sum_{k=1}^{m^2} I_m(\hat{\tau}_k; \tau_k)$$
$$= \sum_{k=1}^{m^2} H_m(\tau_k) - H_m(\tau_k | \hat{\tau}_k) \geq \frac{15}{16} m^2.$$

The second inequality above is due to data processing inequality applied to the fact $\hat{\tau}_k$ only depends on $\hat{\theta}^*$. Finally the claim can be concluded by data processing inequality $I_m(\hat{\theta}; \theta) \geq I_m(\hat{\theta}^*; \theta)$.

CLAIM 2. *Let $\delta > 0$ be a parameter that will be specified later. For any $\delta > 0$, there exist a constant $C_3 > 0$ (depending on $C, \alpha_1$ and $\delta$) such that*

$$\mathbb{P}_m\left(\frac{dQ_m}{dQ_0}(X_1) > C_3\right) \leq \delta, \tag{14}$$

$$I_m\left(X_1; \theta \middle| \frac{dQ_m}{dQ_0}(X_1) > C_3\right) \mathbb{P}_m\left(\frac{dQ_m}{dQ_0}(X_1) > C_3\right) \leq \frac{u}{2m}. \tag{15}$$

PROOF OF CLAIM 2: We first prove (14) holds with large enough constant $C_3$. Let $X_{1,k}$ denote the $k$-th coordinate of $X_1$. Note that

$$\frac{dQ_m}{dQ_0}(X_1) = \prod_{k=1}^{u} \left(e^{-\frac{\lambda^2}{2m}} \cdot \frac{e^{-\frac{\lambda X_{1,k}}{\sqrt{m}\sigma}} + e^{\frac{\lambda X_{1,k}}{\sqrt{m}\sigma}}}{2}\right).$$

Using the basic inequality $\ln(\frac{e^t + e^{-t}}{2}) \leq \frac{t^2}{2}$, we have

$$\mathbb{P}_m\left(\frac{dQ_m}{dQ_0}(X_1) > C_3\right) = \mathbb{P}_m\left(\ln\frac{dQ_m}{dQ_0}(X_1) > \ln C_3\right)$$

$$= \mathbb{P}_m\left(\sum_{k=1}^{u}\left(\ln\left(\frac{e^{-\frac{\lambda X_{1,k}}{\sqrt{m}\sigma}} + e^{\frac{\lambda X_{1,k}}{\sqrt{m}\sigma}}}{2}\right) - \frac{\lambda^2}{2m}\right) > \ln C_3\right)$$

$$\leq \mathbb{P}_m\left(\sum_{k=1}^{u} \frac{\lambda^2}{2m\sigma^2}\left(X_{1,k}^2 - \sigma^2\right) > \ln C_3\right)$$

$$= \mathbb{P}_m\left(\sum_{k=1}^{u} \frac{\lambda^2}{2m\sigma^2}\left(X_{1,k}^2 - \sigma^2 - \frac{\lambda^2\sigma^2}{m}\right) > \ln C_3 - \frac{\lambda^4 u}{2m^2}\right).$$

Note that $\sum_{k=1}^{u} \frac{\lambda^2}{2m\sigma^2}\left(X_{1,k}^2 - \sigma^2 - \frac{\lambda^2\sigma^2}{m}\right)$ has mean 0 and variance at most $(1+\lambda^2)\lambda^4 u/m^2$. Note that we have assumed $M_1 \leq C\sigma m^{2\alpha_1 + \frac{1}{2}}$, this implies $u \leq C^{\frac{2}{2\alpha_1+1}} m^2$. So by Chebyshev's inequality, as long as

$$\ln C_3 \geq \frac{C^{\frac{2}{2\alpha_1+1}}\lambda^4}{2} + \sqrt{C^{\frac{2}{2\alpha_1+1}}(1+\lambda^2)\lambda^4/\delta},$$

we have

$$\mathbb{P}_m\left(\frac{dQ_m}{dQ_0}(X_1) > C_3\right) \leq \delta.$$

We now prove the second inequality (15). Note that when $\theta \sim \mathbb{U}(S_{m,u})$, $\theta$ and event $\{\frac{dQ_m}{dQ_0}(X_1) > C_3\}$ are independent (due to symmetry of $S_{m,u}$). Define

$$\theta_a = (\frac{\lambda\sigma}{\sqrt{m}}, \frac{\lambda\sigma}{\sqrt{m}}, ..., \frac{\lambda\sigma}{\sqrt{m}}, 0, 0, ..., 0) \in S_m.$$

By symmetry, it is easy to show that

$$I_m\left(X_1; \theta \middle| \frac{dQ_m}{dQ_0}(X_1) > C_3\right) \mathbb{P}_m\left(\frac{dQ_m}{dQ_0}(X_1) > C_3\right) = \int_{\frac{dQ_m}{dQ_0}(X_1) > C_3} p(x_1|\theta = \theta_a) \log \frac{p(x_1|\theta = \theta_a)}{q_m(x_1)} dx_1$$

where $p(x_1|\theta = \theta_a)$ denote the density of $x_1$ when $\theta = \theta_0$, and $q_m(x_1)$ denote the density of law $Q_m$.

Further, note that we have following decomposition for $p(x_1|\theta = \theta_a)$ and $q_m(x_1)$:

$$p(x_1|\theta = \theta_a) = \prod_{i=1}^{u} \frac{1}{\sqrt{2\pi}\sigma} e^{-\frac{(x_{1,k} - \frac{\lambda\sigma}{\sqrt{m}})^2}{2\sigma^2}},$$

$$q_m(x_1|\theta = \theta_a) = \prod_{i=1}^{u} \frac{1}{2\sqrt{2\pi}\sigma} \left(e^{-\frac{(x_{1,k} - \frac{\lambda\sigma}{\sqrt{m}})^2}{2\sigma^2}} + e^{-\frac{(x_{1,k} + \frac{\lambda\sigma}{\sqrt{m}})^2}{2\sigma^2}}\right).$$

So we can get
(16)
$$\int_{\frac{dQ_m}{dQ_0}(X_1) > C_3} p(x_1|\theta = \theta_a) \log \frac{p(x_1|\theta = \theta_a)}{q_m(x_1)} dx_1$$

$$= u \int \frac{1}{\sqrt{2\pi}\sigma} e^{-\frac{(y - \frac{\lambda\sigma}{\sqrt{m}})^2}{2\sigma^2}} \cdot \log\left(\frac{2}{1 + \exp(-\frac{2\lambda y}{\sqrt{m}\sigma})}\right) \mathbb{P}_m\left(\frac{dQ_m}{dQ_0}(X_1) > C_3 \middle| x_{1,1} = y\right) dy$$

$$\leq u \int_{y \in [-2\lambda\sqrt{m}\sigma, 2\lambda\sqrt{m}\sigma]} \frac{1}{\sqrt{2\pi}\sigma} e^{-\frac{(y - \frac{\lambda\sigma}{\sqrt{m}})^2}{2\sigma^2}} \cdot \log\left(\frac{2}{1 + \exp(-\frac{2\lambda y}{\sqrt{m}\sigma})}\right) \mathbb{P}_m\left(\frac{dQ_m}{dQ_0}(X_1) > C_3 \middle| x_{1,1} = y\right) dy$$

$$+ u \int_{y \notin [-2\lambda\sqrt{m}\sigma, 2\lambda\sqrt{m}\sigma]} \frac{1}{\sqrt{2\pi}\sigma} e^{-\frac{(y - \frac{\lambda\sigma}{\sqrt{m}})^2}{2\sigma^2}} \cdot \log\left(\frac{2}{1 + \exp(-\frac{2\lambda y}{\sqrt{m}\sigma})}\right) \mathbb{P}_m\left(\frac{dQ_m}{dQ_0}(X_1) > C_3 \middle| x_{1,1} = y\right) dy.$$

Now we bound the first term of the right hand side in (16). It can be shown that when $C_3$ is a large enough constant, we could get

$$\mathbb{P}_m\left(\frac{dQ_m}{dQ_0}(X_1) > C_3 \middle| x_{1,1} = 2\lambda\sqrt{m}\sigma\right) \leq \frac{\ln 2}{4\lambda^2},$$

(we omit the proof here because it is similar to the proof of (14).)

Thus it is easy to show

$$\int_{y\in[-2\lambda\sqrt{m}\sigma, 2\lambda\sqrt{m}\sigma]} \frac{1}{\sqrt{2\pi}\sigma} e^{-\frac{(y-\frac{\lambda\sigma}{\sqrt{m}})^2}{2\sigma^2}} \cdot \log\left(\frac{2}{1+\exp(-\frac{2\lambda y}{\sqrt{m}\sigma})}\right) \mathbb{P}_m\left(\frac{dQ_m}{dQ_0}(X_1) > C_3 \bigg| x_{1,1} = y\right) dy$$

$$\leq \frac{\ln 2}{4\lambda^2} \int \frac{1}{\sqrt{2\pi}\sigma} e^{-\frac{(y-\frac{\lambda\sigma}{\sqrt{m}})^2}{2\sigma^2}} \cdot \log\left(\frac{2}{1+\exp(-\frac{2\lambda y}{\sqrt{m}\sigma})}\right) dy$$

$$\leq \frac{\ln 2}{4\lambda^2} \cdot \frac{1}{\ln 2} \left(\frac{\lambda}{\sqrt{m}}\right)^2 = \frac{1}{4m}.$$

where the second inequality is due to the entropy bound given in Michalowicz et al. (2008).

Next we are going to bound the second term of the right hand side in (16). Because $\lambda \geq 10$, it is easy to show

$$\int_{y\notin[-2\lambda\sqrt{m}\sigma, 2\lambda\sqrt{m}\sigma]} \frac{1}{\sqrt{2\pi}\sigma} e^{-\frac{(y-\frac{\lambda\sigma}{\sqrt{m}})^2}{2\sigma^2}} \cdot \log\left(\frac{2}{1+\exp(-\frac{2\lambda y}{\sqrt{m}\sigma})}\right) \mathbb{P}_m\left(\frac{dQ_m}{dQ_0}(X_1) > C_3 \bigg| x_{1,1} = y\right) dy$$

$$\leq \log 2 \cdot \int_{y\notin[-2\lambda\sqrt{m}\sigma, 2\lambda\sqrt{m}\sigma]} \frac{1}{\sqrt{2\pi}\sigma} e^{-\frac{(y-\frac{\lambda\sigma}{\sqrt{m}})^2}{2\sigma^2}} dy \leq 2\exp\left(-\frac{(2\lambda\sqrt{m} - \frac{\lambda}{\sqrt{m}})^2}{2}\right) < \frac{1}{4m}.$$

Apply the above two bounds to (16) we can get

$$\int_{\frac{dQ_m}{dQ_0}(X_1) > C_3} p(x_1|\theta = \theta_a) \log \frac{p(x_1|\theta = \theta_a)}{q_m(x_1)} dx_1 \leq \frac{u}{2m}$$

when $C_3$ is a large enough constant. This directly implies inequality (15).

Denote the set $R = \{x \in \mathbb{R}^\infty : \frac{dQ_m}{dQ_0}(x) \leq C_3\}$ and random variable $W_i = \mathbb{I}_{\{X_i \in R\}}$.

CLAIM 3. *For each $i = 1, 2, ..., m$, we have*

$$I_m(Z_i; \theta | X_i \in R)\mathbb{P}_m(X_i \in R) \leq \frac{256\lambda^2}{m} \left(\mathbb{E}_m(L_i \mathbb{I}_{\{X_i \in R\}}) + H_m(W_i)\right).$$

PROOF OF CLAIM 3: Let $\tilde{Z}_i$ defined as

$$\tilde{Z}_i \triangleq \begin{cases} Z_i & \text{if } X \in R \\ \star & \text{if } X \notin R \end{cases}$$

where $\star$ is a unique symbol which is different with any 0-1 string.

The following lemma is instrumental to establishing later results.

LEMMA 3 (Multidimensional strong data processing inequality). *Suppose $T = (T^{(1)}, T^{(2)}, ..., T^{(d)})$ be a collection of random variables where each entry is an i.i.d Bernoulli random variable with mean $\frac{1}{2}$. Let $\mu_0$ be a d-dimensional vector and $\Delta > 0$ be a positive real number. Let $X$ be a d-dimensional Gaussian random variable where $X^{(1)}, X^{(2)}, ..., X^{(d)}$ are independent with distribution*

$$X^{(k)} \sim N(\mu_0^{(k)} + T^{(k)}\Delta, \sigma^2).$$

*Let $Z$ be a discrete random variable such that $T \to X \to Z$ is a Markov chain, i.e. $Z \perp T | X$. Then the following multidimensional strong data processing inequality holds:*

$$I(T; Z) \le 64 \left(\frac{\Delta}{\sigma}\right)^2 I(X; Z). \tag{17}$$

Lemma 3 has been proved in Cai and Wei (2020a). For sake of completeness, we provide its proof in the present supplementary material.

Apply Lemma 3 on Markov chain $\theta \to X_i \to \tilde{Z}_i$ where $\theta \sim \mathbb{U}(S_m)$, we have

$$I_m(\theta; \tilde{Z}_i) \le \frac{256\lambda^2}{m} I_m(\tilde{Z}_i; X_i).$$

Note that $W_i \perp \theta$ when $\theta \sim \mathbb{U}(S_{m,u})$, and $W_i$ is determined given $\tilde{Z}$, we have

$$\begin{aligned} I_m(\theta; \tilde{Z}_i) &= I_m(\theta; (\tilde{Z}_i, W_i)) = I_m(\theta; \tilde{Z}_i | W_i) + I_m(\theta; W_i) \\ &= I_m(\theta; \tilde{Z}_i | X_i \in R)\mathbb{P}_m(X_i \in R) + I_m(\theta; \tilde{Z}_i | X_i \notin R)\mathbb{P}_m(X_i \notin R) + I_m(\theta; W_i) \\ &= I_m(\theta; \tilde{Z}_i | X_i \in R)\mathbb{P}_m(X_i \in R). \end{aligned}$$

For similar reasons, we have

$$\begin{aligned} I_m(\tilde{Z}_i; X_i) \le H_m(\tilde{Z}_i) &= H_m(\tilde{Z}_i, W_i) = H_m(\tilde{Z}_i | W_i) + H_m(W_i) \\ &= H_m(\tilde{Z}_i | X_i \in R)\mathbb{P}(X_i \in R) + H_m(\tilde{Z}_i | X_i \notin R)\mathbb{P}(X_i \notin R) + H_m(W_i) \\ &= H_m(Z_i | X_i \in R)\mathbb{P}(X_i \in R) + H_m(W_i) \\ &\le \mathbb{E}_m(L_i | X_i \in R)\mathbb{P}(X_i \in R) + H_m(W_i) \end{aligned}$$

where the latter inequality is due to Shannon's source coding theorem (Shannon, 1948).

Combining the above three formulas yields the desired inequality.

PROOF OF THE MAIN THEOREM:

Note that the region $R$ is "symmetric" where $x \in R$ is equivalent to $|x| \in R$ ($|x|$ is entry-wise absolute value). So $\mathbb{P}(X \in R | \theta)$ is invariant for all $\theta \in S_m$, therefore $W_i \perp \theta$

when $\theta \sim \mathbb{U}(S_m)$. Based on this, for each $i = 1, 2, ..., m$ we have

$$
\begin{aligned}
I_m(Z_i; \theta) &\leq I_m((Z_i, W); \theta) = I_m(Z_i; \theta|W) + I_m(W; \theta) \\
&= I_m(Z_i; \theta|W) \\
&= I_m(Z_i; \theta|X_i \in R)\mathbb{P}_m(X_i \in R) + I_m(Z_i; \theta|X_i \notin R)\mathbb{P}_m(X_i \notin R) \\
&\leq I_m(Z_i; \theta|X_i \in R)\mathbb{P}_m(X_i \in R) + I_m(X_i; \theta|X_i \notin R)\mathbb{P}_m(X_i \notin R) \\
&\leq \frac{256\lambda^2}{m}\left(\mathbb{E}_m(L_i \mathbb{I}_{\{X_i \in R\}}) + H_m(W_i)\right) + \frac{u}{2m}
\end{aligned}
\tag{18}
$$

where the second inequality is due to data processing inequality and the last inequality is derived from Claim 2 and Claim 3.

Taking summation over (18), we have

$$
\frac{256\lambda^2}{m}\left(mH_m(W_1) + \sum_{i=1}^{m}\mathbb{E}_m(L_i\mathbb{I}_{\{X_i \in R\}})\right) + \frac{u}{2} \geq \sum_{i=1}^{m} I_m(Z_i; \theta) \geq I_m(\hat{\theta}; \theta) \geq \frac{15}{16}u
\tag{19}
$$

where the last inequality is due to Claim 1.

Note that for each $i = 1, 2, ..., m$, we have

$$
\mathbb{E}_m(L_i\mathbb{I}_{\{X_i \in R\}}) = \mathbb{E}_0(L_i\mathbb{I}_{\{X_i \in R\}}\frac{dQ_m}{dQ_0}(X_i)) \leq C_3 \mathbb{E}_0(L_i\mathbb{I}_{\{X_i \in R\}}) \leq C_3\mathbb{E}_0(L_i).
$$

Substitute the above inequality into (19) we can get

$$
\mathbb{E}_0(L) = \sum_{i=1}^{m}\mathbb{E}_0(L_i) \geq \frac{1}{C_3}\left(\frac{7}{4096\lambda^2}mu - mH_m(W_1)\right).
$$

Note that $H_m(W_1) \leq -\delta\log\delta - (1-\delta)\log(1-\delta)$. We can always set $\delta$ to a sufficient small constant so that $H_m(W_1) \leq \frac{7}{2048\lambda^2}$. Note that $u \geq 1$, then we can conclude that

$$
\mathbb{E}_0(L) \geq \frac{7}{8192C_3\lambda^2}mu.
$$

Finally, for any $\alpha, p, M > 0$, given the fact that $(0, 0, 0, ..., 0, ...) \in \mathcal{B}_{p,q}^{\alpha}(M)$, we have

$$
\sup_{\theta \in \mathcal{B}_{p,q}^{\alpha}(M)} \mathbb{E}_{\theta}L \geq \mathbb{E}_0(L) \geq \frac{7}{8192C_3\lambda^2}mu \gtrsim \left(\frac{M_1}{\sigma}\right)^{\frac{2}{2\alpha_1+1}}m^{\frac{2\alpha_1+2}{2\alpha_1+1}}.
$$

6.3. *Proof of Theorem 5.* This theorem can be viewed as an extension of Theorem 4. Note that there exists $(\alpha_0, M_0, p_0, q_0) \in \tilde{S}_0$ such that

$$
M_0 = \sigma m^{2\alpha_0 + \frac{1}{2}}.
\tag{20}
$$

Note that for any $\hat{\theta} \in \mathcal{G}(\tilde{S}, C(\cdot))$ and $(\alpha, M, p, q) \in \tilde{S}$, we have

$$\sup_{\theta \in \mathcal{B}_{p,q}^{\alpha}(M)} \mathbb{E}\|\hat{\theta} - \theta\|^2 \leq \tilde{C}(\alpha) M^{\frac{2}{2\alpha+1}} \left(\frac{\sigma^2}{m}\right)^{\frac{2\alpha}{2\alpha+1}},$$

$$\sup_{\theta \in \mathcal{B}_{p_0,q_0}^{\alpha_0}(M)} \mathbb{E}\|\hat{\theta} - \theta\|^2 \leq \tilde{C}(\alpha_0) M_0^{\frac{2}{2\alpha_0+1}} \left(\frac{\sigma^2}{m}\right)^{\frac{2\alpha_0}{2\alpha_0+1}}.$$

Based on above two inequalities and (20), apply Theorem 4, then apply (20) again, we can conclude

$$\sup_{\theta \in \mathcal{B}_{p,q}^{\alpha}(M)} \mathbb{E}L(\hat{\theta}) \gtrsim \left(\frac{M_0}{\sigma}\right)^{\frac{2}{2\alpha_0+1}} m^{\frac{2\alpha_0+2}{2\alpha_0+1}} + \left(\frac{M}{\sigma}\right)^{\frac{2}{2\alpha+1}} m^{\frac{2\alpha+2}{2\alpha+1}} = m^3 + \left(\frac{M}{\sigma}\right)^{\frac{2}{2\alpha+1}} m^{\frac{2\alpha+2}{2\alpha+1}}.$$

## SUPPLEMENTARY MATERIAL

**Supplement A: Supplement to "Distributed Nonparametric Function Estimation: Optimal Rate of Convergence and Cost of Adaptation"**
(doi: url to be specified). In this supplementary material, we prove Theorems 1,2,5, Proposition 1 and the technical lemmas.

*arXiv preprint arXiv:1803.01302.*

DEPARTMENT OF STATISTICS  
THE WHARTON SCHOOL  
UNIVERSITY OF PENNSYLVANIA  
PHILADELPHIA, PENNSYLVANIA 19104  
USA  
E-mail: tcai@wharton.upenn.edu  
      hongjiw@wharton.upenn.edu  
URL: http://www-stat.wharton.upenn.edu/~tcai/